\documentclass[a4paper]{mcom-l16-pre}

\usepackage{KITcolors}  
\usepackage[pdftex]{graphicx}
\usepackage{latexsym}
\usepackage{amsmath,amssymb,amsthm,thmtools}
\usepackage{pgfplots,tikz}
\usetikzlibrary{quotes}
\usetikzlibrary{arrows.meta} 
\usetikzlibrary{patterns}
\usetikzlibrary{spy}
\usetikzlibrary{calc}
\pgfplotsset{compat=1.16} 
\usepackage[algo2e,ruled,vlined,linesnumbered]{algorithm2e}
\usepackage{listings}
\usepackage{dsfont}
\usepackage{float}
\usepackage{xcolor}
\usepackage{bbm}
\usepackage{mathtools} 
\usepackage{amsmath}
\usepackage{subcaption}
\usepackage[shortlabels]{enumitem}
\usepackage{hyperref}
\hypersetup{
	colorlinks,
	linkcolor={red!50!black},
	citecolor={blue!50!black},
	urlcolor={blue!70!black}
}
\usepackage[nameinlink]{cleveref}
\usepackage{standalone}
\usepackage[msc-links]{amsrefs}

\crefdefaultlabelformat{\normalfont #2#1#3}
\crefformat{equation}{\normalfont \textup{(#2#1#3)}}
\crefrangeformat{equation}{\normalfont (#3#1#4) to~(#5#2#6)}
\crefmultiformat{equation}{\normalfont (#2#1#3)}{\normalfont\ and~(#2#1#3)}{\normalfont , (#1#1#3)}{\normalfont\ and~(#2#1#3)}

\usepackage{ifthen}
\usepackage{xspace}

\makeatletter
\DeclareFontEncoding{LS2}{}{\noaccents@}
\makeatother
\DeclareFontSubstitution{LS2}{stix2}{m}{n}
\DeclareSymbolFont{largesymbolsSTIX}{LS2}{stix2ex}{m}{n}
\DeclareMathDelimiter{\lParen} {\mathopen} {largesymbolsSTIX}{"DE}{largesymbolsSTIX}{"02}
\DeclareMathDelimiter{\rParen} {\mathclose}{largesymbolsSTIX}{"DF}{largesymbolsSTIX}{"03}
\DeclareMathDelimiter{\lBrack} {\mathopen} {largesymbolsSTIX}{"E0}{largesymbolsSTIX}{"06}
\DeclareMathDelimiter{\rBrack} {\mathclose}{largesymbolsSTIX}{"E1}{largesymbolsSTIX}{"07}
\DeclareMathDelimiter{\lBrace} {\mathopen} {largesymbolsSTIX}{"E8}{largesymbolsSTIX}{"0E}
\DeclareMathDelimiter{\rBrace} {\mathclose}{largesymbolsSTIX}{"E9}{largesymbolsSTIX}{"0F}
\DeclareMathDelimiter{\lbrbrak}{\mathopen} {largesymbolsSTIX}{"EE}{largesymbolsSTIX}{"14}
\DeclareMathDelimiter{\rbrbrak}{\mathclose}{largesymbolsSTIX}{"EF}{largesymbolsSTIX}{"15}
\DeclareSymbolFont{stmry}{U}{stmry}{m}{n}
\DeclareMathDelimiter\llbracket{\mathopen}{stmry}{"4A}{stmry}{"71}
\DeclareMathDelimiter\rrbracket{\mathclose}{stmry}{"4B}{stmry}{"79}

\newlength\harmoscheight
\newlength\harmoscwidth
\newcommand{\abs}[1]{\left|#1\right|}
\newcommand{\MLnorm}[1]{\norm{#1}} %
\newcommand{\bigMLnorm}[1]{\bignorm{#1}} %
\newcommand{\BigMLnorm}[1]{\Bignorm{#1}} %
\newcommand{\norm}[1]{\lVert #1 \rVert}

\newcommand{\bignorm}[1]{\big\lVert #1 \big\rVert}
\newcommand{\Bignorm}[1]{\Big\lVert #1 \Big\rVert}

\newcommand{\mfrac}[2]{\scalebox{0.8}{$\dfrac{#1}{#2}$}{}}
\newcommand{\pmat}[1]{\begin{pmatrix} #1 \end{pmatrix}}
\newcommand{\pmatT}[1]{\begin{pmatrix} #1 \end{pmatrix}^\top}

\newcommand{\dist}{\mathrm{dist}}
\newcommand{\weight}{\omega}

\newcommand{\wz}{\weight_{dz}}
\newcommand{\Iwz}{\NodalInt \wz}
\newcommand{\Ipd}{d_h}
\newcommand{\IpdwzE}{\varepsilon_{dz}}

\newcommand{\CN}{Crank--Nicolson\xspace}
\newcommand{\LF}{leapfrog\xspace}
\newcommand{\DS}{Domain splitting\xspace}
\newcommand{\ds}{domain splitting\xspace}

\newcommand{\timestep}{time step\xspace}
\newcommand{\timesteps}{time steps\xspace}
\newcommand{\timestepsize}{time-step size\xspace}

\newcommand{\R}{\mathbb{R}} %
\newcommand{\N}{\mathbb{N}} %
\newcommand{\dintx}{\, \mathrm{d} x}

\newcommand{\dint}[1]{\, \mathrm{d} #1}
\newcommand{\e}{\mathrm{e}}

\renewcommand{\d}[1]{\partial_{#1}}
\newcommand{\dt}{\d{t}}

\newcommand*{\divg}{\nabla\cdot}

\newcommand{\grad}{\nabla}

\newcommand{\Rminus}{R_{-}}
\newcommand{\Rplus}{R_{+}}
\newcommand{\RCN}{R}
\newcommand{\RminusLF}{\widehat{R}_{-}}
\newcommand{\RplusLF}{\widehat{R}_{+}}
\newcommand{\RLF}{\widehat{R}}

\newcommand{\T}{T}      %
\newcommand{\Ti}{(0,\T]}%
\newcommand{\ov}{\delta}%
\newcommand{\ovp}{\ell} %
\newcommand{\ts}[1][]{%
	\ifthenelse{ \equal{#1}{} }
	{\ensuremath{\tau}}
	{\ensuremath{\frac{\tau}{#1}}}
}

\newcommand{\tss}[1][]{%
	\ifthenelse{ \equal{#1}{} }
	{\ensuremath{\tau^2}}
	{\ensuremath{\frac{\tau^2}{#1}}}
}
\newcommand{\numberTS}{n_T} %

\newcommand{\fulldomain}{\Omega}
\newcommand{\dfulldomain}{\partial \fulldomain}
\newcommand{\fulldomainc}{\overline{\fulldomain}}
\newcommand{\D}{\Theta}

\newcommand{\SD}[1]{\fulldomain_{#1}}
\newcommand{\dSD}[1]{\partial \SD{#1}}
\newcommand{\SDc}[1]{\overline{\fulldomain}_{#1}}

\newcommand{\SDov}[1]{\fulldomain_{#1}^{\ov}}
\newcommand{\dSDov}[1]{\partial \SDov{#1}}
\newcommand{\SDovc}[1]{\overline{\fulldomain}_{#1}^{\ov}}

\newcommand{\numberSD}{N}

\newcommand{\SDovint}[1]{\Gamma_{#1}^{\ov}}
\newcommand{\SDint}[1]{\Gamma_{#1}}
\newcommand{\Avgint}{\Gamma}

\newcommand{\trialfunc}[1][]{
	\ifthenelse{ \equal{#1}{1} }
	{q} %
	{p} %
}
\newcommand{\trialfunch}[1][]{
	\trialfunc[#1]_h
}

\newcommand{\xcoord}{c} %
\newcommand{\solutionVector}{x}  %
\newcommand{\qComponent}{u} %
\newcommand{\pComponent}{v} %

\newcommand{\solVecExact}[1][]{\solutionVector(t_{#1})}
\newcommand{\qExakt}[1][]{\qComponent^{#1}}
\newcommand{\pExakt}[1][]{\pComponent^{#1}}

\newcommand{\solDiscr}{\solutionVector_h} 
\newcommand{\qDiscr}{\qComponent_h}
\newcommand{\pDiscr}{\pComponent_h}

\newcommand{\solDS}[1]{\solutionVector_{\normalfont\text{DS}}^{#1}} %
\newcommand{\qDS}[1]{\qComponent_{\normalfont\text{DS}}^{#1}}
\newcommand{\pDS}[1]{\pComponent_{\normalfont\text{DS}}^{#1}}

\newcommand{\solDSloc}[2]{\solutionVector^{#2}_{#1}} %
\newcommand{\qDSloc}[2]{\qComponent_{#1}^{#2}}
\newcommand{\pDSloc}[2]{\pComponent_{#1}^{#2}}

\newcommand{\solCN}[1]{\solutionVector_{\normalfont\text{CN}}^{#1}}
\newcommand{\qCN}[1]{\qComponent_{\normalfont\text{CN}}^{#1}}
\newcommand{\pCN}[1]{\pComponent_{\normalfont\text{CN}}^{#1}}

\newcommand{\solCNmod}[1]{\widetilde{\solutionVector}_{\normalfont\text{CN}}^{\,#1}}
\newcommand{\qCNmod}[1]{\widetilde{\qComponent}_{\normalfont\text{CN}}^{\,#1}}                %
\newcommand{\pCNmod}[1]{\widetilde{\pComponent}_{\normalfont\text{CN}}^{\,#1}}                %

\newcommand{\qLF}[1]{\normalfont \qComponent_{\text{lf}}^{#1}}
\newcommand{\pLF}[1]{\normalfont \pComponent_{\text{lf}}^{#1}}
\newcommand{\solLF}[1]{\solutionVector_{\normalfont\text{lf}}^{#1}}

\newcommand{\solLFmod}[1]{\widetilde{\solutionVector}_{\normalfont\text{lf}}^{\,#1}}
\newcommand{\qLFmod}[1]{\widetilde{\qComponent}_{\normalfont\text{lf}}^{\,#1}}                %
\newcommand{\pLFmod}[1]{\widetilde{\pComponent}_{\normalfont\text{lf}}^{\,#1}}                %

\newcommand{\qDiffDSCN}{z_{\qComponent,i}^{n}}

\newcommand{\pDiffDSCN}{z_{\pComponent,i}^{n}}

\newcommand{\fmh}[1]{\overline{f}_h^{#1 -1/2}}
\newcommand{\gmh}[1]{\overline{g}_h^{#1 -1/2}}

\newcommand{\Vspace}{V}
\newcommand{\Hspace}{H}
\newcommand{\indexNormV}[1][]{
	\ifthenelse{ \equal{#1}{} }
	{1}
	{1,#1}
}

\newcommand{\indexNormH}[1][]{
	\ifthenelse{ \equal{#1}{} }
	{\fulldomain}
	{#1}
}

\newcommand{\Vh}[1][]{
	\ifthenelse{ \equal{#1}{} }
	{\Vspace_h}
	{\Vspace_h(#1)}
}

\newcommand{\Hh}[1][]{
	\ifthenelse{ \equal{#1}{} }
	{\Hspace_h}
	{\Hspace_h(#1)}
}

\newcommand{\Vhn}[1][]{
	\ifthenelse{ \equal{#1}{} }
	{\Vspace_{h,0}}
	{\Vspace_{h,0}(#1)}
}

\newcommand{\Hhn}[1][]{
	\ifthenelse{ \equal{#1}{} }
	{\Hspace_{h,0}}
	{\Hspace_{h,0}(#1)}
}

\newcommand{\NormVh}[1][]{
	\ifthenelse{ \equal{#1}{} }
	{\Vh}
	{\Vh[#1]}
}

\newcommand{\indexNormVh}[1][]{
	\ifthenelse{ \equal{#1}{} }
	{\Vh}
	{\Vh[#1]}
}

\newcommand*{\Vhnorm}[2][]{\norm{#2}_{\indexNormVh[#1]}}
\newcommand*{\bigVhnorm}[2][]{\bignorm{#2}_{\indexNormVh[#1]}}
\newcommand*{\BigVhnorm}[2][]{\Bignorm{#2}_{\indexNormVh[#1]}}

\newcommand{\NormHh}[1][]{
	\ifthenelse{ \equal{#1}{} }
	{\Hh}
	{\Hh[#1]}
}
\newcommand{\indexNormHh}[1][]{
	\ifthenelse{ \equal{#1}{} }
	{\Hh}
	{\Hh[#1]}
}

\newcommand*{\Hhnorm}[2][]{\MLnorm{#2}_{\indexNormHh[#1]}}
\newcommand*{\bigHhnorm}[2][]{\bigMLnorm{#2}_{\indexNormHh[#1]}}
\newcommand*{\BigHhnorm}[2][]{\BigMLnorm{#2}_{\indexNormHh[#1]}}

\newcommand{\Espace}[1][]{
	\ifthenelse{ \equal{#1}{} }
	{X}
	{X(#1)}
}

\newcommand{\dEspace}[1][]{
	\ifthenelse{ \equal{#1}{} }
	{X_h}
	{X_h(#1)}
}

\newcommand{\indexNormEspace}[1][]{
	\ifthenelse{ \equal{#1}{} }
	{\Espace}
	{\Espace[#1]}
}

\newcommand{\indexNormdEspace}[1][]{
	\ifthenelse{ \equal{#1}{} }
	{\dEspace}
	{\dEspace[#1]}
}

\newcommand*{\Espacenorm}[2][]{\norm{#2}_{\indexNormEspace[#1]}}
\newcommand*{\bigEspacenorm}[2][]{\bignorm{#2}_{\indexNormEspace[#1]}}

\newcommand*{\dEspacenorm}[2][]{\norm{#2}_{\indexNormdEspace[#1]}}
\newcommand*{\bigdEspacenorm}[2][]{\bignorm{#2}_{\indexNormdEspace[#1]}}
\newcommand*{\BigdEspacenorm}[2][]{\Bignorm{#2}_{\indexNormdEspace[#1]}}

\newcommand{\projdEspace}{\mathcal{J}_h}

\newcommand{\errorDSCN}[1]{E^{#1}}

\newcommand{\Cavg}{C_{\normalfont\text{avg}}}
\newcommand{\cavg}{C_{\normalfont\text{avg},\Hh}}
\newcommand{\CML}{C_{\normalfont\text{ML}}}
\newcommand{\Cdata}{C_{\normalfont\text{data}}}  %
\newcommand{\Cglob}{C_{\normalfont\text{glob}}}
\newcommand{\MaxCFLellc}{M_\ovp}
\newcommand{\MaxCFLellcs}{\MaxCFLellc^2}
\newcommand{\Cint}{C_{\!\NodalInt}}
\newcommand{\CInt}{C_{\normalfont\text{Int}}}
\newcommand{\CRitz}{C_{\Rproj}}
\newcommand{\Cinv}{C_{\normalfont\text{inv}}}
\newcommand{\Derror}{D_{\ts,h,\ov}}

\newcommand{\calculatingstep}{n}
\newcommand{\mLapl}{L}  %
\newcommand{\dLapl}{L_{h}}  %
\newcommand{\dsysMat}{A_{h}} %
\newcommand{\id}{\mathrm{I}}

\newcommand{\averaging}{\zeta}
\newcommand{\vectoraveraging}{\zeta} %
\newcommand{\avOpDS}{\vectoraveraging \Big( \big\{ \solDSloc{i}{\calculatingstep} \big\}_{i=1,\dots, \numberSD} \Big)}

\newcommand{\avOpDiff}{\vectoraveraging \Big( \big\{  \solDSloc{i}{\calculatingstep} - \restr{\solCNmod{\calculatingstep}}{\SDov{i}} \big\}_{i=1,\dots, \numberSD} \Big)}

\def\L#1{{L^{#1}(\fulldomain)}}
\def\LD#1{{L^{#1}(\D)}}

\def\H#1{{H^{#1}(\fulldomain)}}
\def\Hn#1{{H^{#1}_{0}(\fulldomain)}}

\newcommand{\Th}[1][]{
	\ifthenelse{ \equal{#1}{} }
	{\mathcal{T}_h}
	{\mathcal{T}_h(#1)}
}
\newcommand{\K}{K}

\newcommand{\nodalPoints}[1][]{
	\ifthenelse{ \equal{#1}{} }
	{\mathcal{N}}
	{\mathcal{N}_{#1}}
}

\newcommand{\diam}{\mathrm{diam}}

\newcommand{\FemS}[1]{W_h(#1)}
\newcommand{\FemSn}[1]{W_{h,0}(#1)}
\newcommand{\FemSfull}{\FemS{\fulldomain}}

\newcommand{\FemSnfull}{\FemSn{\fulldomain}}

\newcommand{\Rproj}{\mathcal{R}_h}
\newcommand{\NodalInt}{\mathcal{I}_h}

\newcommand{\ps}{\kappa}

\newcommand{\ip}[1]{\big( #1 \big)}

\newcommand{\ipD}[2]{\left({#2}\right)_{#1}}

\newcommand{\ipSDov}[2]{\ipD{\SDov{#1}}{#2}}
\newcommand{\ipfD}[1]{\ipD{\fulldomain}{#1}}

\newcommand{\ipML}[2]{ \big\lParen#2\big\rParen_{#1}}

\newcommand{\ba}[1]{a_{#1}}
\newcommand{\baSD}[1]{\ba{\SD{#1}}}
\newcommand{\baSDov}[1]{\ba{\SDov{#1}}}

\newcommand{\baw}[2]{\ba{#1,#2}}

\newcommand{\bb}[1]{b_{#1}}

\newcommand{\bbSDov}[1]{\bb{\SDov{#1}}}

\newcommand{\bbw}[2]{\bb{#1,#2}}

\newcommand{\sumSD}{\sum_{i=1}^{\numberSD}}

\newcommand{\sumTh}[1][]{
	\ifthenelse{ \equal{#1}{} }
	{\sum_{K \in \Th}}
	{\sum_{K \in \Th[#1]}}
}

\renewcommand{\dim}{m}

\newcommand\restr[2]{{%
\left.\kern-\nulldelimiterspace %
#1 %
\vphantom{\big|} %
\right|_{#2} %
}}

\newcommand{\interior}[1]{%
	{\kern0pt#1}^{\mathrm{o}}%
}

\newcommand{\landauO}[1]{\mathcal{O}(#1)}
\newcommand{\half}{1/2}

\newcommand{\card}{\mathop{\mathrm{card}}}
\newcommand{\CSI}{Cauchy--Schwarz inequality\xspace}

\newtheorem{theorem}{Theorem}[section]
\newtheorem{lemma}[theorem]{Lemma}

\newtheorem{definition}[theorem]{Definition}

\ifpdf
\hypersetup{
  pdftitle={A domain splitting method for the acoustic wave equation},
  pdfauthor={T. Buchholz, M. Hochbruck}
}
\fi

\begin{document}
\title[A domain splitting method for the acoustic wave equation]{A non-iterative domain decomposition time integrator for linear wave equations}

\author[Tim Buchholz]{Tim Buchholz}
\address{Institute for Applied and Numerical Mathematics,
    Karlsruhe Institute of Technology, Englerstr. 2, 76131 Karls\-ru\-he, Germany}
\email{\{tim.buchholz,marlis.hochbruck\}@kit.edu}
\thanks{Funded by the Deutsche Forschungsgemeinschaft (DFG, German Research Foundation) -- Project-ID 258734477 -- SFB 1173}

\author[Marlis Hochbruck]{Marlis Hochbruck}


\begin{abstract}
	We propose and analyze a non-iterative domain decomposition integrator for the linear acoustic wave equation. 
	The core idea is to combine an implicit \CN step on spatial subdomains with a local prediction step at the subdomain interfaces. 
	This enables parallelization across space while advancing sequentially in time, without requiring iterations at each time step. 
	The method is similar to the methods in \cites{BluLR92,DawD92}, which have been designed for parabolic problems.
	Our approach adapts them to the case of the wave equation in a fully discrete setting, using linear finite elements with mass lumping. 
	Compared to explicit schemes, our method permits significantly larger time steps and retains high accuracy. 
	We prove that the resulting method achieves second-order accuracy in time and global convergence of order $\mathcal{O}(h + \tau^2)$ under a CFL-type condition, which depends on the overlap width between subdomains. 
	We conclude with numerical experiments which confirm the theoretical results.
\end{abstract}

\maketitle

\section{Introduction} 
\label{sec::introduction}
Wave equations are central to modeling a wide range of physical phenomena, including acoustic, electromagnetic, and elastic wave propagation, cf. \cite{DoeHKRSW23}*{Chapter 1}. 
Numerical simulations of such problems pose significant challenges, particularly due to the computational cost associated with resolving wave dynamics accurately over large spatial and temporal domains.

A critical aspect of simulating wave propagation is the time integration process. %
Broadly, time-stepping methods fall into two categories: explicit and implicit methods. 
Explicit methods are only conditionally stable:
Step size restrictions, often also called Courant-Friedrichs-Lewy (CFL) conditions, force the choice of the time step to directly depend (linearly or worse) on the space discretization width. 
Unfortunately, the need for numerous tiny steps can lead to a significant computational burden, slowing down the simulation process.
Implicit methods, in contrast, are unconditionally stable in many cases and permit larger time steps.
However, implicit methods incur a higher computational cost per step, as they require solving a system of equations at each step.
This cost becomes particularly pronounced when dealing with three-dimensional or finely resolved two-dimensional problems. 

To address these computational challenges, we propose a non-iterative domain decomposition (DD) method for time integration of a finite element discretization of the linear acoustic wave equation.
Specifically, we consider the first-order formulation of the linear acoustic wave equation on an open, polygonal domain $\fulldomain \subset \mathbb{R}^\dim$, given by
    \begin{align} 
        \label{eq::waveEqFirstOrder}
        \begin{cases}
            \begin{array}{rlrlll}
                \dt \qExakt                            & = \pExakt ,            & \dt \pExakt  & = \mLapl u  +  f \quad & \text{in }  \fulldomain \times \Ti, \\
                \qExakt(x,0)                           & = \qExakt[0](x), \quad & \pExakt(x,0) & = \pExakt[0](x)                          & \text{in }  \fulldomain,            \\
                \restr{\qExakt(\cdot,t)}{\dfulldomain} & = 0                    &              &                                          & \text{for } t \in  \Ti,             \\
            \end{array}
        \end{cases}
    \end{align}
    with $\mLapl u = \divg(\ps^2 \grad u) $,
    a constant wave propagation speed $\ps > 0$ and final time $T$.
	The differential operator $\mLapl$ is defined on its domain $ D(\mLapl)= \H2 \cap \Hn1 $.
For the space discretization, we consider linear finite elements with mass lumping on a triangular mesh.

Domain decomposition methods are well-established tools for solving spatially discretized partial differential equations (PDEs) by splitting the spatial domain into subdomains, often enabling parallel computation. 
Classical Schwarz methods iteratively solve stationary problems by exchanging boundary data between subdomains.  
An overview of such methods can be found in \cite{GanH13} and \cite{GanZ22}, respectively.

For time-dependent problems, Schwarz waveform relaxation (SWR) methods extend this concept by integrating the PDE independently in each subdomain and exchanging interface data iteratively \cites{GanHN99,GanHN03,HalS09,GanKM21}. 
However, their iterative nature introduces a computational overhead that conflicts with the goal of reducing the computational cost of the time integration.

We follow a non-iterative approach motivated by \cite{BluLR92} for parabolic problems which 
is based on an overlapping decomposition of the spatial domain. 
Coupling between the subdomains is established via artificial boundary conditions obtained from local extrapolation at the interfaces in each time step.
The authors of \cite{BluLR92} show that this \textit{\ds} method is stable and exhibits optimal order of convergence if the overlap regions between the subdomains are sufficiently large. 
Related non-iterative methods are discussed in \cites{DawD92,DawD94}.
Further work on time-integrators based on domain decomposition 
can be found under the names \textit{domain decomposition operator splitting} or \textit{regionally-additive schemes}; see for instance \cites{Vab08,VabZa13,EfeVa21} and \cites{HanEr17,EisHa18}, respectively. Most of these works focus on the parabolic case. Recently the authors of \cite{BucDjEi24} managed to extend such a method to stochastic evolution equations.

Another non-iterative DD approach for hyperbolic problems is tent pitching \cites{DraGSW22,GopMS15,GopSW17,CiaGM24}
in which space-time domains are constructed incrementally with (typically tent-shaped) subdomains, adhering to causality constraints. 
Yet, this can be complex for non-constant coefficients and lacks straightforward parallelization due to subdomain solution dependencies.
A different approach relying on finite propagation speed and superposition principles was recently proposed in  \cite{GalMa23}. 
Similar to our convergence analysis, the authors of  \cite{GalMa23} show convergence against the \CN scheme and utilize  exponential decay in space.
However, their localized implicit time stepping method requires significantly larger overlap regions than the scheme presented here.

In this work, we adapt the \ds method from \cite{BluLR92} to the setting of the linear acoustic wave equation in first order formulation. 
In this method, which we call \textit{\ds} method hereafter, we divide the spatial domain into overlapping subdomains on which the \CN method is applied independently. 
We use a local variant of the \LF scheme for the prediction of boundary data on the subdomain interfaces.

We prove that the proposed \ds method achieves second-order accuracy in time and overall convergence of order $\landauO{h + \ts^2}$ under a CFL-type condition that depends on the wave propagation speed, the spatial mesh size, and the width of the overlap between subdomains.
Compared to the explicit \LF scheme, the proposed method allows $\landauO{\ell}$  bigger time steps, 
where $\ell$ describes the number of mesh layers within the overlap.
At the same time the \ds method retains convergence of order $\landauO{h + \ts^2}$, given sufficient regularity of the exact solution, the initial data, and the source term.

\bigskip

The paper is organized as follows.
In \Cref{sec::discretization} we first introduce some notation and discuss standard space and second-order time discretizations on the full spatial domain. 
We then construct the \ds scheme in \Cref{sec::DomainSplitting} and present the main result in \Cref{sec::MainResult}. 
The subsequent error analysis spans \Cref{sec::PredictionError,sec::LocalError,sec::AvgProps,sec::ErrorRecursion}.
In \Cref{sec::PredictionError}, we bound the error introduced by the predictions at the boundaries of the subdomains.
In \Cref{sec::LocalError}, we show that the local error can be interpreted as a solution of a stationary problem, from which we deduce an exponential decay in space.  
\Cref{sec::AvgProps} gathers some properties of the averaging process, which is necessary to obtain a globally continuous approximation. 
In \Cref{sec::ErrorRecursion}, we derive an error recursion and conclude the error analysis with the proof of the main result.
Finally, we present numerical experiments in \Cref{sec::NumericalExp}.

\section{Discretization and notation}
\label{sec::discretization}

\subsection{Discretization in space}
Let $ \Th = \Th[\fulldomain] $ be a (shape- and contact-regular) matching simplicial mesh of $ \fulldomain $, see, e.g. \cite{ErnG21a}*{Definition 8.11}. The parameter $h$ denotes the minimal diameter of the elements in $\Th$.
We use continuous linear finite elements and define the approximation spaces
\begin{align*}
	\FemSfull   = \{ \trialfunc[1] \in \H1 : \forall \K \in \Th: \restr{\trialfunc[1]}{\K} \in P_1 \} \quad \text{and} \quad 
	\FemSnfull  = \FemSfull \cap \Hn1,
\end{align*}
where $P_1 $ denotes the set of all linear polynomials in $\dim$ variables.

Let $\D \subseteq \fulldomain$ such that
\begin{equation}
	\label{eq::defAdmissibleSubdomain}
	\K \cap \D = \K \quad \text{or} \quad  \K\cap \D = \emptyset \quad \text{for all } K\in \Th.
\end{equation}
We denote the submesh containing all elements in $\D$ by $\Th[\D]$. 

Moreover, we define the set of all nodal points of the triangulation $ \Th $
by
\begin{equation*}
	\nodalPoints[\fulldomainc] \coloneqq  \bigcup\limits_{K \in \Th} \{\xcoord_{K,j}\}_{j=1}^{\dim + 1} ,
\end{equation*}
where $ \{\xcoord_{K,j}\}_{j=1}^{\dim+1} $ are the vertices of the element $ K \in \Th $.
The set of interior nodal points of $\D \subseteq \fulldomain$ satisfying \eqref{eq::defAdmissibleSubdomain} is given by
$\nodalPoints[\D] = \nodalPoints[\fulldomainc] \cap \D $ and for $\D = \fulldomain$ we write $\nodalPoints = \nodalPoints[\fulldomain]$.

\subsection{Norms and bilinear forms}
\label{sssec::NormsBilinearForms}

For $ \D \subseteq \fulldomain $ with \eqref{eq::defAdmissibleSubdomain} we denote the standard $ \LD2 $ inner product by $ \ipD{\D}{\cdot,\cdot} = \ipD{\LD2}{\cdot,\cdot} $ and define the bilinear form 
\begin{equation}
	\label{def::ba}
	\ba{\D}\ip{\trialfunc[1],\trialfunc[2]} \coloneqq \int_{\D} \ps^2 \grad \trialfunc[1] \cdot \grad \trialfunc[2] \dintx, \qquad \text{for } \trialfunc[1],\trialfunc[2] \in H^1(\D).
\end{equation}
We formulate \eqref{eq::waveEqFirstOrder} weakly by seeking $\pmatT{\qComponent, \pComponent} \in \Hn1 \times \L2$ such that
\begin{equation}
    \label{eq::waveEqFirstOrder-weakly}
    \begin{aligned}
        \ipfD{\dt \qComponent, \varphi} &= \ipfD{ \pComponent, \varphi}, &&\text{for all } \varphi \in H_0^1(\fulldomain) ,\\
    \ipfD{ \dt \pComponent, \psi} &= \ba{\fulldomain} \ip{u, \psi}  + \ipfD{f, \psi}, && \text{for all } \psi \in H_0^1(\fulldomain) .
    \end{aligned}
\end{equation}
Moreover, we use mass-lumping, i.e., we approximate the $L^2$ inner product cell-wise by a (Lobatto-)quadrature approximation, see, e.g., \cite{CohPe17}.
If we denote the nodal interpolation on a single cell $K$ by $ \NodalInt^K $, we can express the mass-lumped inner product by
\begin{equation*}
	\ipML{\D}{\trialfunch[1],\trialfunch[2]} = \sumTh[\D] \int_K \NodalInt^K (\trialfunch[1] \trialfunch[2]) \dintx , \qquad \trialfunch[1],\trialfunch[2] \in \FemS{\D}.
\end{equation*}
On a mesh consisting of intervals, triangles, or tetrahedrons we can give the mass-lumped bilinear form $ \ipML{\D}{\cdot,\cdot} $ explicitly by the following representation
\begin{equation}
    \label{eq::massLumping-explicitRepr}
	\ipML{\D}{\trialfunch[1],\trialfunch[2]} = \sum_{K\in \Th[\D]} \sum_{j=1}^{\dim+1} \frac{\abs{K}}{\dim+1} \trialfunch[1](\xcoord_{K,j}) \trialfunch[2](\xcoord_{K,j}),
\end{equation}
see, e.g., \cite{Cia02}*{Chapter 4}.
The mass-lumped inner product induces a norm on $ \FemSfull $, which is 
equivalent to the standard $ \L2 $-norm, i.e.
\begin{equation}
	\label{const::CML}
	\CML^{-1} \ipML{\fulldomain}{\trialfunch[1],\trialfunch[1]}^{\half} \leq \norm{\trialfunch[1]}_{\L2} \leq \CML \ipML{\fulldomain}{\trialfunch[1],\trialfunch[1]}^{\half}, \qquad \trialfunch[1] \in \FemSfull,
\end{equation}
see, e.g., \cite{Tho97}*{Chapter 15} and \cite{Rav73}*{Section 3 and 4}. %

We then introduce the finite-dimensional Hilbert spaces  $\Vh[\D]$, $\Hh[\D]$ by 
\begin{equation*}
	\Vh[\D] = \bigl( \FemS{\D}, \ba{\D}\ip{\cdot, \cdot} \bigr), \qquad \Hh[\D] = \bigl( \FemS{\D}, \ipML{\D}{\cdot, \cdot} \bigr)
\end{equation*}
with norms 
\begin{equation*}
	\Vhnorm[\D]{\cdot}^2 = \ba{\D}\ip{\cdot, \cdot} \quad\text{and}\quad \Hhnorm[\D]{\cdot}^2 = \ipML{\D}{\cdot, \cdot},
\end{equation*}
respectively. 
The subspaces with a zero trace are denoted as $\Vhn[\D]$ and $\Hhn[\D]$.
The linear operator $ \dLapl: \Vhn[\fulldomain] \to \Vhn[\fulldomain] $ associated to $ \ba{\fulldomain}\ip{\cdot, \cdot} $ is given by
\begin{equation}
	\label{def::Lh}
	\ipML{\fulldomain}{\dLapl \trialfunch[1],\trialfunch[2]} = \ba{\fulldomain}\ip{\trialfunch[1],\trialfunch[2]}, \quad \trialfunch[1],\trialfunch[2] \in \Vhn[\fulldomain]. 
\end{equation}

For $ \Espace[\D] = H^1(\D)\times\LD2$ and $ \solutionVector = \pmatT{\qComponent , \pComponent} \in \Espace[\D] $ we define
\begin{equation*}
	\bigEspacenorm[\D]{\solutionVector}^2 =  \ba{\D}\ip{\qComponent,\qComponent}  +\ipD{\D}{\pComponent,\pComponent}.
\end{equation*}
Similarly, for $ \dEspace[\D] =  \Vh[\D] \times \Hh[\D] $ and $ \solDiscr = \pmatT{\qDiscr, \pDiscr} \in \dEspace[\D] $ we denote the discrete variant using the mass-lumped scalar product by
\begin{equation*}
	\bigdEspacenorm[\D]{\solDiscr}^2  = \Vhnorm[\D]{\qDiscr}^2  + \Hhnorm[\D]{\pDiscr}^2 .
\end{equation*}
If $ \D =  \fulldomain $,  we drop the set $ \D $ and write, e.g.,
$
	\Espacenorm{\cdot} = \Espacenorm[\fulldomain]{\cdot} $ and $
	\dEspacenorm{\cdot} = \dEspacenorm[\fulldomain]{\cdot}.
$
We use the $L^2$-projection with respect to $\ipML{\fulldomain}{\cdot,\cdot}$ to approximate the right-hand side 
$f_h(t) \approx f(t)$. 
We want to emphasize, that this projection coincides with the nodal interpolation $\NodalInt$, as 
\begin{equation}
	\label{eq::MLprojectIsInterpolation}
	\ipML{\fulldomain}{\trialfunc[1] - \NodalInt \trialfunc[1], w_h} = 0, \quad\text{for }  \trialfunc[1] \in C(\fulldomain), w_h \in \Hh.
\end{equation}
For $\trialfunc[1] \in H^2(K)$ on $K \in \Th$, the standard interpolation estimate
\begin{equation}
	\label{const::CIh}
	\abs{\trialfunc[1] - \NodalInt \trialfunc[1]}_{1,K} \leq \Cint h \abs{\trialfunc[1]}_{2,K}
\end{equation}
holds with $\abs{\cdot}_{\alpha,K}$ denoting the standard $H^{\alpha}$-seminorms on $K$, see, e.g. \cite{BreS08}*{Theorem 4.4.4}.
We also use the nodal interpolation $\NodalInt$ for the initial values, so that we have
\begin{equation}  \label{eq:projected-data}
	\qComponent_h^0 = \NodalInt \qExakt[0], \qquad \pComponent_h^0 = \NodalInt \pExakt[0], \qquad f_h^n = \NodalInt f(t_n).
\end{equation}
In $\dEspace$ the discrete counterpart of $\eqref{eq::waveEqFirstOrder-weakly}$ becomes finally
\begin{equation}
	\label{eq::wave-semidiscrete}
	\dt \pmat{\qComponent_h \\ \pComponent_h} = \pmat{0 & I \\ - \dLapl & 0} \pmat{\qComponent_h \\ \pComponent_h} + \pmat{0 \\ f_h(t)}.
\end{equation}
Note, that the discrete problem is stated in terms of the discrete operator $\dLapl$ on the space $\dEspace[]$. For the implementation, we need to assemble mass and stiffness matrices from a basis of  $\Vh$ and $\Hh$.

\subsection{Classical second-order time integration}
\label{subsec::clasicalTimeIntegrators}
Let $\ts>0$ be a given \timestepsize and $t_n = n \ts $ for $n= 0,1,\dots,\numberTS$ with final time $T =\numberTS\ts$. 
The implicit \CN scheme for the semi-discrete system \eqref{eq::wave-semidiscrete} yields approximations $\solCN{n} = \pmatT{\qCN{n} , \pCN{n}}\in \dEspace$ given by 
\begin{subequations}
	\label{FullDiscrCN-ML}
	\begin{alignat}{1}
		\label{FullDiscrCN-ML-a}
		\qCN{n} &= \qCN{n-1} + \ts[2]\big( \pCN{n} + \pCN{n-1} \big) ,\\
		\label{FullDiscrCN-ML-b}
		\pCN{n} &= \pCN{n-1} - \ts[2]\dLapl \big( \qCN{n} + \qCN{n-1} \big) + \ts\fmh{n},
	\end{alignat}
\end{subequations}
where we denote $\overline{f}_h^{n-1/2} = \frac12(f_h^n + f_h^{n-1})$. For the implementation it is advantageous to solve 
\begin{equation}
	\label{eq::ReformulationCN_qComp}
	\big( I + \tss[4] \dLapl \big) \qCN{n} = \big( I - \tss[4] \dLapl \big) \qCN{n-1} + \ts \pCN{n-1} + \tss[2]\fmh{n}
\end{equation}
for $\qCN{n}$ instead of solving the coupled system \eqref{FullDiscrCN-ML}.

The one-step formulation of the explicit leapfrog scheme applied to \eqref{eq::wave-semidiscrete} yields an approximation 
$\solLF{n} = \pmatT{\qLF{n} , \pLF{n}}\in \dEspace$
given by 
\begin{subequations}
	\label{FullDiscrLF-ML}
	\begin{alignat}{1}
		\qLF{n - \half} &= \qLF{n-1} + \ts[2] \pLF{n-1} \label{FullDiscrLF-ML-a}\\
		\pLF{n} &= \pLF{n-1} - \ts[] \dLapl \qLF{n -\half} + \ts \fmh{n} \label{FullDiscrLF-ML-b}\\
		\qLF{n} &= \qLF{n-\half} + \ts[2] \pLF{n} \label{FullDiscrLF-ML-c}
	\end{alignat}
\end{subequations}
If we choose a nodal basis of $\FemSfull$ and assemble the mass and stiffness matrices for \eqref{FullDiscrLF-ML}, we get one linear system with a diagonal mass matrix in each step. We would like to stress that this makes the method local in space.
It is well known that the leapfrog scheme suffers from a strong CFL condition, namely 
\begin{equation}
	\label{eq::CFL-leapfrog}
	\tss \bignorm{\dLapl}_{\Hh \leftarrow \Hh} \leq 4,
\end{equation}
see, e.g., \cite{Jol03}*{Theorem 1}.

\section{\DS scheme}  
\label{sec::DomainSplitting}

The construction of our new scheme is inspired by the non-iterative  DD scheme from
\cite{BluLR92} proposed for finite element discretizations of
parabolic problems. 
We start by decomposing the spatial domain into non-overlapping subdomains $\SD{i}$,
$i=1,\ldots,\numberSD$, such that
\begin{equation*}
	\overline{\fulldomain} = \bigcup_{i=1}^N \overline{\SD{i}}.
\end{equation*}

Next, we define an overlapping decomposition with overlapping subdomains $\SDov{i}$ by extending $\SD{i}$ with $\ovp$ layers of elements. Here, $\ov$ denotes the minimal width of the layers extending $\SD{i}$, see \Cref{fig:A4_decomposition}. 
Obviously, we have $\ov \sim \ovp h$ for a regular mesh.

Solving the discrete wave equation \eqref{eq::wave-semidiscrete} on each subdomain $\SDov{i}$ requires appropriate boundary conditions. Hence, one \timestep of the \ds method consists of the following three steps:
\begin{enumerate}[1)]
	\item %
		\textbf{Prediction:}
	      We generate boundary values at all artificial boundaries $\dSDov{i} \cap \fulldomain$ by
		  a leapfrog step, which is local in space.
	\item %
		\textbf{\CN on subdomains:}
	      Using the predicted boundary values, we perform one \timestep of the \CN{} method \eqref{FullDiscrCN-ML} on each of the overlapping subdomains $\SDov{i}$.
	\item%
		\textbf{Averaging:}
	      In this step, we restrict the subdomain solution on $\SDov{i}$ to the non--overlapping subdomain $\SD{i}$, $i=1,\ldots,N$. On the interfaces $ \dSD{i}\cap \dSD{j} $, we apply averaging to obtain a unique continuous solution on $\fulldomain$.  
\end{enumerate} 

\begin{figure}[t]
	\centering
	\begin{subfigure}[t]{0.5\textwidth}
        \centering
        \begin{tikzpicture}[scale=0.26] 
			\def\maxX{20}
			\def\maxY{20}
			\def\ELL{3}

			\definecolor{colorSD1}{RGB}{0,0,103}
			\definecolor{colorSD2}{RGB}{0,184,208}
			\definecolor{colorSD3}{RGB}{226,199,0}
			\definecolor{colorSD4}{RGB}{113,0,0}

			\draw[draw=colorSD1, fill=colorSD1,opacity = 0.5] (0,0) rectangle (\maxX*0.5,\maxY*0.5);
			\draw[draw=colorSD2, fill=colorSD2,opacity = 0.5] (0,\maxY*0.5) rectangle (\maxX*0.5,\maxY);
			\draw[draw=colorSD3, fill=colorSD3,opacity = 0.5] (\maxX*0.5,0) rectangle (\maxX,\maxY*0.5);
			\draw[draw=colorSD4, fill=colorSD4,opacity = 0.5] (\maxX*0.5,\maxY*0.5) rectangle (\maxX,\maxY);
			\def\opa{0.25}
		
			\draw[line width=1pt] (0,0) rectangle (\maxX,\maxY);
			\foreach \x in {0,1,...,\maxX}{
					\draw[opacity=\opa] (\x,0) -- (\x,\maxY);
				}
			\foreach \y in {0,1,...,\maxY}{
					\draw[opacity=\opa] (0,\y) -- (\maxX,\y);
				}
			\foreach \x in {1,...,\maxX}{
					\foreach \y in {1,...,\maxY}
					\draw[opacity=\opa] (\x,\y) -- (\x-1,\y-1);
				}
			\draw[line width=2pt,opacity=\opa] (\maxX*0.5,0) -- (\maxX*0.5,\maxY);
			\draw[line width=2pt,opacity=\opa] (0,\maxY*0.5) -- (\maxX,\maxY*0.5);
			\draw[draw=colorSD1, line width=0.5, fill=colorSD1, opacity = 0.4] (0,0) rectangle (\maxX*0.5+\ELL,\maxY*0.5+\ELL);
			\node[circle,draw=colorSD1,text=black,fill=white,fill opacity=0.75,inner sep = 0.5mm,scale=1.] at (\maxX*0.25+\ELL*0.5,\maxY*0.25+\ELL*0.5) {$ \Omega_1^{\delta} $};
			\node[circle,draw=colorSD1,text=black,fill=white,fill opacity=0.75,inner sep = 0.5mm,scale=1.] at (\maxX*0.75,\maxY*0.25) {$ \Omega_2 $};
			\node[circle,draw=colorSD1,text=black,fill=white,fill opacity=0.75,inner sep = 0.5mm,scale=1.] at (\maxX*0.25,\maxY*0.75) {$ \Omega_3 $};
			\node[circle,draw=colorSD1,text=black,fill=white,fill opacity=0.75,inner sep = 0.5mm,scale=1.] at (\maxX*0.75,\maxY*0.75) {$ \Omega_4 $};
			\draw[color=colorSD1,line width = 2.5] (0.05,\maxY*0.5+\ELL) -- (\maxX*0.5+\ELL,\maxY*0.5+\ELL);
			\draw[color=colorSD1,line width = 2.5] (\maxX*0.5+\ELL,0.05) -- (\maxX*0.5+\ELL,\maxY*0.5+\ELL);
			\node[circle,draw=white,text=black,fill=white,fill opacity=0.75,inner sep = 0.5mm,scale=1.] at (\maxX*0.25+\ELL*0.5+0.5,\maxY*0.5+\ELL*0.5) {$ \delta $};
			\draw[<->,color=black!30!white,line width = 1.] (\maxX*0.25+0.5,\maxY*0.5) -- (\maxX*0.25+0.5,\maxY*0.5+\ELL);
		\end{tikzpicture}
    \end{subfigure}%
	\begin{subfigure}[t]{0.5\textwidth}
        \centering
        \begin{tikzpicture}[scale=0.75]

\coordinate (P0) at (0.0,0.0);
\coordinate (P1) at (8.0,0.0);
\coordinate (P2) at (8.0,5.0);
\coordinate (P3) at (0.0,5.0);
\coordinate (P4) at (0.7272727272713526,0.0);
\coordinate (P5) at (1.454545454542336,0.0);
\coordinate (P6) at (2.181818181812941,0.0);
\coordinate (P7) at (2.909090909083548,0.0);
\coordinate (P8) at (3.636363636354153,0.0);
\coordinate (P9) at (4.363636363626778,0.0);
\coordinate (P10) at (5.090909090901422,0.0);
\coordinate (P11) at (5.818181818176066,0.0);
\coordinate (P12) at (6.545454545450711,0.0);
\coordinate (P13) at (7.272727272725355,0.0);
\coordinate (P14) at (8.0,0.7142857142849608);
\coordinate (P15) at (8.0,1.428571428568591);
\coordinate (P16) at (8.0,2.142857142852804);
\coordinate (P17) at (8.0,2.857142857139607);
\coordinate (P18) at (8.0,3.571428571426405);
\coordinate (P19) at (8.0,4.285714285713203);
\coordinate (P20) at (7.272727272724246,5.0);
\coordinate (P21) at (6.545454545453539,5.0);
\coordinate (P22) at (5.818181818185859,5.0);
\coordinate (P23) at (5.09090909091818,5.0);
\coordinate (P24) at (4.3636363636505,5.0);
\coordinate (P25) at (3.636363636378783,5.0);
\coordinate (P26) at (2.909090909103026,5.0);
\coordinate (P27) at (2.181818181827269,5.0);
\coordinate (P28) at (1.454545454551513,5.0);
\coordinate (P29) at (0.7272727272757562,5.0);
\coordinate (P30) at (0.0,4.285714285714285);
\coordinate (P31) at (0.0,3.57142857142857);
\coordinate (P32) at (0.0,2.857142857142853);
\coordinate (P33) at (0.0,2.142857142857135);
\coordinate (P34) at (0.0,1.428571428571423);
\coordinate (P35) at (0.0,0.7142857142857117);
\coordinate (P36) at (0.6406971084116695,2.496786780131428);
\coordinate (P37) at (7.360322976148707,2.495589387098663);
\coordinate (P38) at (4.726161661606117,4.369587365236362);
\coordinate (P39) at (3.274726537514178,0.6287207136905271);
\coordinate (P40) at (2.544133927692886,4.370283672790174);
\coordinate (P41) at (5.454545454538743,0.6298366572994336);
\coordinate (P42) at (6.174140979179779,4.370689754642936);
\coordinate (P43) at (1.818181818177639,0.6298366572959363);
\coordinate (P44) at (6.867470873281542,0.6086417541664666);
\coordinate (P45) at (1.090909090913635,4.370163342699604);
\coordinate (P46) at (7.3494028615226,3.855735598708358);
\coordinate (P47) at (0.6547688164076052,1.139225694939238);
\coordinate (P48) at (1.806623297078301,4.369135070301849);
\coordinate (P49) at (2.165552901591474,3.74401452210516);
\coordinate (P50) at (2.903259317152124,3.74106067674268);
\coordinate (P51) at (2.531041954038659,3.118460716884476);
\coordinate (P52) at (3.264539607384775,3.112664202679856);
\coordinate (P53) at (2.900225990673177,2.489008913070915);
\coordinate (P54) at (3.62755704109661,2.48470550851066);
\coordinate (P55) at (2.167993293669921,2.492724346204462);
\coordinate (P56) at (3.991815802230595,3.110589057540085);
\coordinate (P57) at (4.353459784589963,2.484667396614063);
\coordinate (P58) at (4.717943890312323,3.109181614686769);
\coordinate (P59) at (5.079683227766391,2.479851588050829);
\coordinate (P60) at (5.443626727814732,3.106498214282779);
\coordinate (P61) at (5.80486699124863,2.478568030807522);
\coordinate (P62) at (4.716804475017599,1.862641115030466);
\coordinate (P63) at (2.539899902798354,1.86943739441172);
\coordinate (P64) at (6.178358572172582,3.106254332226279);
\coordinate (P65) at (1.812861992943986,1.867350784092799);
\coordinate (P66) at (6.555890540127068,2.482492574583282);
\coordinate (P67) at (6.173082117646791,1.857608065862945);
\coordinate (P68) at (6.887810002541,1.851267606102768);
\coordinate (P69) at (1.421631961522508,2.496508572865697);
\coordinate (P70) at (1.078321678387604,1.838228385733559);
\coordinate (P71) at (1.04518871783308,3.11692851310739);
\coordinate (P72) at (6.919242762810557,3.147726354150095);
\coordinate (P73) at (6.550056927047473,3.767467815863466);
\coordinate (P74) at (5.813892855468084,3.742435862118718);
\coordinate (P75) at (1.782257816415151,3.12387245431642);
\coordinate (P76) at (4.357101524241143,3.738564152152282);
\coordinate (P77) at (3.270487493896649,4.368748608779173);
\coordinate (P78) at (6.172024206958466,0.6199998102283633);
\coordinate (P79) at (5.45098581289999,4.369780147215225);
\coordinate (P80) at (5.447590848832858,1.858767811924908);
\coordinate (P81) at (2.546397198701698,0.6264179262653663);
\coordinate (P82) at (4.724112902303668,0.6242213325309243);
\coordinate (P83) at (3.996646733988828,0.6267715268962997);
\coordinate (P84) at (3.631405401422447,1.250869768280416);
\coordinate (P85) at (5.084952078723731,3.739341225948689);
\coordinate (P86) at (3.997221699792107,4.369583040829473);
\coordinate (P87) at (3.630862861924443,3.738521225859338);
\coordinate (P88) at (1.44452847870528,1.2215096578411);
\coordinate (P89) at (2.909396085873609,1.250536363762852);
\coordinate (P90) at (2.174027437031106,1.24657683510206);
\coordinate (P91) at (6.522374425257523,1.236656264708224);
\coordinate (P92) at (3.263178703252664,1.874386774249048);
\coordinate (P93) at (5.815040261984466,1.238261741191983);
\coordinate (P94) at (5.083749536477915,1.243663313938261);
\coordinate (P95) at (3.989208881955642,1.870588460054185);
\coordinate (P96) at (1.416582876566147,3.746484001266361);
\coordinate (P97) at (0.6504769182568116,3.824834745188991);
\coordinate (P98) at (4.353392356431248,1.246971245035119);
\coordinate (P99) at (6.894596821138778,4.412210941267588);
\coordinate (P100) at (1.103916492366339,0.5851999219475496);
\coordinate (P101) at (7.320836105287467,1.165185466531415);
\coordinate (P102) at (7.51063692170502,1.807842800367948);
\coordinate (P103) at (0.4747575206413757,1.809133886446557);
\coordinate (P104) at (7.524013956604464,3.214285714283007);
\coordinate (P105) at (0.4661104944371756,3.21428571428571);
\coordinate (P106) at (7.465452751953227,0.5311612121781986);
\coordinate (P107) at (0.5345472480476442,4.468838787821019);
\coordinate (P108) at (7.475798340905039,4.479372478390775);
\coordinate (P109) at (0.524201659093822,0.5206275216090231);

\definecolor{colorSD1}{RGB}{0,0,103}

\fill[draw=colorSD1, fill=colorSD1,opacity = 0.8] (P33) -- (P36) -- (P53) -- (P92) -- (P84) -- (P7) -- (P0) -- cycle;

\fill[draw=colorSD1, fill=colorSD1,opacity = 0.5] (P33) -- (P32) -- (P31) -- (P97) -- (P96) -- (P49) -- (P50) -- (P87) -- (P56) -- (P57) -- (P62) -- (P94) -- (P82) -- (P9) -- (P8) -- (P7) -- (P39) -- (P84) -- (P92) -- (P53) -- (P55) -- (P69) -- (P36) -- cycle;

\foreach \a/\b/\c in {
            P72/P46/P73,P14/P15/P101,P70/P47/P88,P14/P101/P106,P46/P72/P104,P37/P66/P68,P47/P70/P103,P66/P37/P72,P69/P36/P70,P36/P69/P71,P73/P46/P99,P30/P31/P97,P88/P47/P100,P49/P48/P96,P55/P51/P75,P64/P66/P72,P64/P72/P73,P75/P49/P96,P69/P55/P75,P51/P49/P75,P37/P68/P102,P30/P97/P107,P73/P42/P74,P71/P69/P75,P48/P45/P96,P64/P73/P74,P71/P96/P97,P74/P42/P79,P81/P89/P90,P74/P79/P85,P76/P38/P86,P79/P38/P85,P12/P44/P78,P70/P36/P103,P40/P48/P49,P38/P76/P85,P25/P26/P77,P42/P22/P79,P40/P27/P48,P22/P23/P79,P82/P94/P98,P38/P24/P86,P64/P61/P66,P65/P55/P69,P65/P69/P70,P59/P62/P80,P60/P61/P64,P28/P29/P45,P12/P13/P44,P43/P6/P81,P28/P45/P48,P23/P38/P79,P66/P67/P68,P60/P59/P61,P41/P11/P78,P25/P77/P86,P63/P55/P65,P5/P6/P43,P21/P22/P42,P10/P41/P82,P26/P40/P77,P11/P12/P78,P53/P51/P55,P76/P86/P87,P26/P27/P40,P10/P11/P41,P86/P77/P87,P40/P49/P50,P24/P25/P86,P66/P61/P67,P27/P28/P48,P39/P83/P84,P7/P39/P81,P6/P7/P81,P58/P59/P60,P9/P10/P82,P39/P8/P83,P50/P49/P51,P58/P57/P59,P9/P82/P83,P61/P59/P80,P59/P57/P62,P53/P55/P63,P8/P9/P83,P7/P8/P39,P23/P24/P38,P60/P64/P74,P58/P60/P85,P56/P57/P58,P52/P53/P54,P89/P63/P90,P56/P54/P57,P52/P54/P56,P56/P58/P76,P50/P51/P52,P52/P51/P53,P76/P58/P85,P40/P50/P77,P60/P74/P85,P94/P62/P98,P50/P52/P87,P77/P50/P87,P67/P61/P80,P56/P76/P87,P52/P56/P87,P18/P19/P46,P34/P35/P47,P16/P17/P37,P32/P33/P36,P42/P73/P99,P81/P39/P89,P65/P70/P88,P78/P44/P91,P93/P80/P94,P43/P81/P90,P83/P82/P98,P53/P63/P92,P65/P88/P90,P62/P57/P95,P80/P62/P94,P67/P80/P93,P88/P43/P90,P39/P84/P89,P63/P65/P90,P62/P95/P98,P92/P84/P95,P63/P89/P92,P41/P93/P94,P20/P21/P99,P4/P5/P100,P68/P67/P91,P71/P75/P96,P54/P53/P92,P41/P78/P93,P82/P41/P94,P89/P84/P92,P54/P92/P95,P57/P54/P95,P101/P15/P102,P84/P83/P98,P96/P45/P97,P95/P84/P98,P72/P37/P104,P91/P67/P93,P43/P88/P100,P78/P91/P93,P36/P71/P105,P91/P44/P101,P21/P42/P99,P5/P43/P100,P13/P1/P106,P29/P3/P107,P1/P14/P106,P3/P30/P107,P68/P91/P101,P16/P37/P102,P34/P47/P103,P36/P33/P103,P18/P46/P104,P37/P17/P104,P32/P36/P105,P2/P20/P108,P0/P4/P109,P19/P2/P108,P35/P0/P109,P101/P44/P106,P71/P97/P105,P97/P45/P107,P15/P16/P102,P33/P34/P103,P97/P31/P105,P17/P18/P104,P31/P32/P105,P4/P100/P109,P20/P99/P108,P68/P101/P102,P44/P13/P106,P45/P29/P107,P46/P19/P108,P47/P35/P109,P99/P46/P108,P100/P47/P109}{
   \draw[thin,opacity=0.5] (\a) -- (\b) -- (\c) -- cycle;
}

\node[circle,draw=colorSD1,text=black,fill=white,fill opacity=0.75,inner sep = 0.5mm,scale=1.] at (P63) {$ \Omega_1^{\delta} $};
\draw[color=colorSD1,line width = 2.0] (P31)  -- (P97)  -- (P96)  -- (P49)  -- (P50)  -- (P87)  -- (P56)  -- (P57)  -- (P62)  -- (P94)  -- (P82)  -- (P9) ; 
\draw[<->, color=black!30!white,line width = 1.] ($(P92)$) -- node[circle,draw=white,text=black,fill=white,fill opacity=0.75,inner sep = 0.5mm,scale=1.,xshift=0.5em,yshift=-1em] {$\delta$} ($(P57)$);

\end{tikzpicture}
    \end{subfigure}%
	\caption{Non-overlapping and overlapping decomposition for $\ovp = 3$ on an equidistant mesh and for $\ovp = 2$ on a non-equidistant mesh.}
	\label{fig:A4_decomposition} 
\end{figure}

We now present the three steps in detail.

\subsection{Prediction}

To generate boundary conditions at each of the interfaces $\SDovint{i} = \dSDov{i} \cap \fulldomain$ we perform a \LF step. This step is fully explicit and only uses the degrees of freedom surrounding the interface $\SDovint{i}$ since mass lumping leads to a diagonal matrix in \eqref{FullDiscrLF-ML}.

In contrast to our approach, in \cite{BluLR92} extrapolation in time was used locally at the desired degrees of freedom.
Using a leapfrog step for the prediction is advantageous here, since it leads to a one-step scheme which simplifies the implementation as well as the analysis.

\subsection{\CN\ step on subdomains}

The prediction step provides inhomogeneous Dirichlet boundary conditions on all interfaces $\SDovint{i}$, $i=1,\ldots,N$. 
This enables us to apply a \CN\ step on each subdomain $\SDov{i}$.
Using the \CN\ method on subdomains is considerably cheaper than on the full domain $\fulldomain$, since the linear systems of equations are smaller and can be solved in parallel.

\subsection{Averaging}
Recall that $ \nodalPoints[\SD{i}] $ is the set of all interior nodes in $ \SD{i} $, i.e., $\nodalPoints[\SD{i}] = \nodalPoints\cap \SD{i}$. We define the averaging operator
\begin{equation*}
	\averaging: \bigtimes\limits_{i = 1}^{\numberSD} \dEspace[\SDov{i}] \to \dEspace, \quad \text{or} \quad  \averaging: \bigtimes\limits_{i = 1}^{\numberSD} \Vh[\SDov{i}] \to \Vh,
\end{equation*} 
mapping a set of subdomain functions $\trialfunc[1]_i\in \dEspace[\SDov{i}]$ or $\trialfunc[1]_i\in \Vh[\SDov{i}]$, $i=1,\dots,\numberSD$, to a global approximation 
by its nodal values
\begin{equation}
		\label{eq::averaging:mean}
		\averaging \bigl( \{\trialfunc[1]_i\}_{i=1,\dots,\numberSD} \bigr)(\xcoord_j) \! = \!  
		\begin{cases} 
			\hspace*{-0.4cm}
			\begin{array}{rll}
				 & \trialfunc[1]_i(\xcoord_j),                                    & \xcoord_j \in \nodalPoints[\SD{i}] \text{ for exactly one } i \in \{1,\dots, \numberSD\},                              \\
				 & \sum\limits_{k \in J_j } \! \frac{\trialfunc[1]_k(\xcoord_j)}{\abs{J_j}},  & \xcoord_j \in \nodalPoints \setminus \big( \mathop{\bigcup}_{i=1,\dots,\numberSD} \nodalPoints[\SD{i}]\big) , \\[-0.5em]
				 &&\; \text{with} \; J_{j} =  \big\{ i \in \{1,\dots, \numberSD \} : \xcoord_j \in \SDc{i}\big\}.
			\end{array}
		\end{cases}
\end{equation}
Note that $\xcoord_j$ is contained in at most one set of interior nodal points $\nodalPoints[\SD{i}]$ as the subdomains $\SD{i}, i = 1, \dots , \numberSD,$ are non-overlapping. 
Moreover, for $\trialfunc[1] \in \Vh[\fulldomain]$ we immediately see
	\begin{equation}
		\label{eq::avg_prop1}
		\averaging \big( \big\{ \restr{\trialfunc[1]}{\SDov{i}}\big\}_{i=1,\dots,\numberSD} \big) = \trialfunc[1],
	\end{equation}
i.e., first taking restrictions to the overlapping subdomains and then applying $ \averaging $ keeps functions in $\Vh[\fulldomain]$ invariant.

\subsection{Full algorithm}
We summarize the \ds method in \Cref{algorithm::DSAlgWave}.
For the initial values we assume $\qExakt[0], \pExakt[0] \in  \H2 $ to give the point-wise evaluations within the interpolation a proper meaning.

\begin{algorithm2e}[ht]
	\DontPrintSemicolon
	\SetAlgoLined
	\SetKwInput{Input}{input}
	\SetKwInOut{Output}{result}
	\vspace{0.5cm}
	\Input{initial data $ \qDS{0}  = \NodalInt \qExakt[0] \in \Vhn$, $ \pDS{0}  = \NodalInt \pExakt[0] \in \Hh $, \\
		overlapping decomposition $ \overline{\fulldomain} = \bigcup\limits_{i=1}^{\numberSD} \SDovc{i} $,\\
		\timestepsize $\ts$, final time $ \T = \numberTS \ts   $} 
	\vspace{0.25cm}
	$ n = 1 $\;
	\While{$n \leq \numberTS$}{
		\textbf{prediction step}: prediction by $ \qLFmod{n}  $ (leapfrog step with i.c. $\solDS{n-1}$) \;
		\For{i = 1,\dots, $ \numberSD $}{
			\textbf{\CN on subdomains:} \\
			calculate $ \solDSloc{i}{n} $ on $ \SDov{i} $ with  $ \restr{\qLFmod{n}}{\SDovint{i}} $ as b.c. on $ \SDovint{i} $\;
		}
		\textbf{averaging step:} construct global $ \solDS{n} = \avOpDS  $ according to \eqref{eq::averaging:mean}\;
		update $ n \leftarrow n+1 \; $
	}
	\Output{$\solDS{\numberTS} \approx \solutionVector(\cdot, T)$.} %
	\vspace{0.5cm}
	\caption{\DS scheme}
	\label{algorithm::DSAlgWave}
\end{algorithm2e} %

\subsection{Further notation}

For the subsequent error analysis, we introduce some additional notation. By $ \solDS{n} $ and $ \solDSloc{i}{n} $ we denote the approximations of \Cref{algorithm::DSAlgWave}. In addition, the \CN\ solution after $n$ steps started from $\solutionVector^0 = \pmat{\qExakt[0], \pExakt[0]}$ is denoted by $ \solCN{n} $. Moreover, we will make use of \textit{local in time} approximations, which are obtained by doing one time step with the \CN\ or the \LF scheme starting from the domain splitting approximation $ \solDS{n-1} $. These approximations are denoted by $ \solCNmod{n} $ and $ \solLFmod{n} $, respectively.

We summarize the notation in \Cref{table::notation:approx}.

\begin{center} 
	\begin{table}[ht]
	\begin{tabular}{||c | c | c | c  | c ||} 
	 \hline
	 domain  	& subdomain &  & \textit{local in time}  & \textit{local in time} \\
	 splitting	&  approx. & \CN & \CN& \LF \\
	        	&  in $n$th step & & &  \\[0.5ex]  
	 \hline\hline
	 $\solDS{n} = \pmat{\qDS{n} \\ \pDS{n}}$ 
	 & $\solDSloc{i}{n} = \pmat{\qDSloc{i}{n} \\ \pDSloc{i}{n}}$ 
	 & $\solCN{n} = \pmat{\qCN{n} \\ \pCN{n}}$
	 & $\solCNmod{n} = \pmat{\qCNmod{n} \\ \pCNmod{n}}$ 
	 & $\solLFmod{n} = \pmat{\qLFmod{n} \\ \pLFmod{n}}$
	 \\ [1ex]
	 \hline
	\end{tabular}
	\caption{Notation for different approximations}
	\label{table::notation:approx}
	\end{table}
\end{center}

\section{Main result}
\label{sec::MainResult}

A main contribution of this work is to show that the \ds method  for approximating the solution of \eqref{eq::wave-semidiscrete} satisfies an error bound of the form
\begin{equation*}
		\bigEspacenorm{\solDS{n} - \solVecExact[n]} \lesssim \ts^2 + h
\end{equation*}
under suitable conditions on the exact solution $\solutionVector $ of \eqref{eq::waveEqFirstOrder}, the right-hand side $f$, and the parameters $h,\ts,$ and $\ovp$.

A key idea is to start with comparing the \ds method and the \CN method, since for the latter, rigorous error bounds of the full discretization exist, see \cite{HipHS19} and \cite{Hip17}*{Example 5.3, Corollary 5.9}.
More precisely, we will show 
\begin{equation*}
    \errorDSCN{n} =	\bigdEspacenorm{\solDS{n} - \solCN{n}} \lesssim \ts^2
\end{equation*}
in the following lemma.

\begin{lemma}%
	\label{ConvergenceDS}
	Let $ \qExakt[0],\pExakt[0] \in  \H2 \cap \Hn1  $ and $ f \in L^{\infty}\bigl([0,T]; \H2 \cap \Hn1\bigr) $ and consider the \ds method of the previous section with a sufficiently large overlap $ \ov \sim \ovp h $.
	Assume that the CFL condition 
	\begin{equation}
		\label{eq::CFLconditionPrediction}
		\tss \bignorm{\dLapl}_{\Hh \leftarrow \Hh } \leq 4 \ovp^2
	\end{equation}
	holds for a \timestep $\ts \in (0,1]$.
	Then there exist constants $ 0 \leq \sigma  \leq 1 $ and $\Cdata>0$ such that
	\begin{align*}
		\errorDSCN{n} & = \bigdEspacenorm{\solDS{n} - \solCN{n}}
		\leq \ts^2  \Cdata\min \{e^{\sigma t_{n}}-1 , \sigma t_n \e^{\sigma \ts t_n} \}
	\end{align*}
	The constant $\Cdata$ is independent of $h$ and $\tau$, for details see \Cref{lem::dLaplonCNsol} below.
\end{lemma}

The proof of this result is quite involved and postponed to \Cref{sec::ErrorAnalysis}.

\begin{theorem}[Error bound against the exact solution]  \label{thm::main}
	Let the assumptions of \Cref{ConvergenceDS} be fulfilled and let the solution of \eqref{eq::waveEqFirstOrder} suffice 
	\begin{equation*}
		u \in C^4\bigl([0,T];\L2\bigr)\cap C^3\bigl([0,T];\Hn1\bigr) \cap C^2\bigl([0,T]; \H2\bigr) \cap C\bigl([0,T];D(\mLapl^2)\bigr).
	\end{equation*} 
	Then, the error of the \ds scheme satisfies
	\begin{align*}
		\bigEspacenorm{\solDS{n} - \solVecExact[n]} \lesssim \ts^2 + h.
	\end{align*}
\end{theorem}  

\begin{proof}
	Let $\projdEspace$ be the projection from $\Espace$ to $\dEspace$ obtained by taking a Ritz projection $\Rproj$ (see, e.g., \cite{ErnG21b}*{Section 32.4}) in both components. Then, by the triangle inequality, we can estimate
	\begin{equation*}
		\bigEspacenorm{\solVecExact[n]- \solDS{n} } \! \leq \! \bigEspacenorm{\solVecExact[n] - \projdEspace \solVecExact[n] } \! +  
		\CML\bigdEspacenorm{\projdEspace\solVecExact[n] - \solCN{n}} \! + 
		\CML \bigdEspacenorm{\solCN{n} - \solDS{n}} \!,
	\end{equation*}
	with $\CML>0$ from \eqref{const::CML}.
	For the first term we can apply standard approximation results to obtain
	\begin{align*}
		\bigEspacenorm{\solVecExact[n] - \projdEspace \solVecExact[n] } \! &= \!\Bigl( \big|{\qComponent(t_n) - \! \Rproj \qExakt(t_n)}\big|_{\H1}^2 \! + \bignorm{\dt \qExakt(t_n) - \! \Rproj (\dt \qExakt(t_n))}_{\L2}^2 \! \Bigr)^{\! \half} \\
		&\lesssim h \Bigl(  \bignorm{\qExakt(t_n)}_{\H2}^2 + \bignorm{\dt \qExakt(t_n)}_{\H1}^2\Bigr)^{\half} , 
	\end{align*}
	 see, e.g., \cite{ErnG21b}*{Theorem 32.15 and 33.2}.
	For the second term, we have
	\begin{align*}
		\bigdEspacenorm{\solCN{n}- \projdEspace\solVecExact[n]} \lesssim& \sqrt{2} \Bigl(  
			 h \norm{\qExakt[0]}_{\H2} +  h \norm{\pExakt[0]}_{\H1}
			+ \tss[8]  \int_{0}^{t_{n+1}} \Espacenorm{\dt^3 \solutionVector(s)} \dint{s} \\
			&+  h \ts \sum_{j=1}^{n} \max\limits_{s \in [t_j,t_{j+1}]} \norm{\dt^2 \qComponent(s)}_{\H1} \\
			&+  h \ts \sum_{j=1}^{n} \norm{\mLapl u (t_j)}_{\H2} + \norm{\mLapl u (t_{j+1})}_{\H2}
		\Bigr), 
	\end{align*}
	from \cite{HipHS19}*{Section 4.8} and \cite{Hip17}*{Section 5.3 and Corollary 5.9}.
 
	Finally, the bound for $\bigdEspacenorm{\solCN{n} - \solDS{n}}$ follows from \Cref{ConvergenceDS}.
\end{proof}

\section{Prediction error}
\label{sec::PredictionError}

A key part of the error analysis is estimating the prediction error at each step, which is challenging because each step begins with a solution that has been cut and reassembled along the interface, offering little regularity. To address this, we aim to analyze the prediction error globally in space, but only locally in time. This is done by using equivalent formulations of the \LF and \CN methods on $\dEspace[]$, %
which can be found in \cite{DoeHKRSW23}*{Sections~11.1 and 11.2}.

\begin{lemma}  \label{lem:LF-CN-onestep}
    \begin{enumerate}[a),leftmargin=\parindent,align=left,labelwidth=\parindent,labelsep=0pt]
        \item The \CN method \eqref{FullDiscrCN-ML} is equivalent to
        \begin{equation}
        	\Rminus \pmat{\qCN{n} \\ \pCN{n}} =
        	\Rplus \pmat{\qCN{n-1} \\ \pCN{n-1}} + \ts \pmat{0 \\ \fmh{n}} \; ,
        	\label{eq::CN_on_VH_operators}
        \end{equation}
    	where
        \begin{equation*}
            \Rminus = \pmat{\id & -\ts[2] \id \\ \ts[2] \dLapl & \id}, \qquad \Rplus = \pmat{\id \quad & \ts[2] \id \\ -\ts[2] \dLapl & \id}.
        \end{equation*}
        \item The \LF method \eqref{FullDiscrLF-ML} is equivalent to
        \begin{equation}
        	\RminusLF \pmat{\qLF{n} \\ \pLF{n}} =
        	\RplusLF \pmat{\qLF{n-1} \\ \pLF{n-1}} 
        	+ \ts \pmat{0 \\ \fmh{n}},
        	\label{eq::LF_on_VH_operators}
        \end{equation}
    		where
        \begin{equation*}
            \RminusLF 
            = \pmat{\id  & -\ts[2] \id \\ \ts[2] \dLapl & \id  - \tss[4] \dLapl} 
            \qquad \RplusLF 
            = \pmat{\id & \ts[2] \id \\ -\ts[2] \dLapl & \id - \tss[4] \dLapl } 
        \end{equation*}
    \end{enumerate}
    \label{def::operatorNotation}
\end{lemma}

Moreover, we define
\begin{equation}
	\label{eq::def_R_operators}
	\RCN = \Rminus^{-1} \Rplus \quad \text{and} \quad \RLF = \RminusLF^{-1} \RplusLF.
\end{equation}	
	
 We start with a stability result for the \CN method, which makes use of the inverse inequality
 \begin{equation}  \label{eq:inverse-ineq}
 	\Hhnorm{\dLapl\trialfunch[1]} \leq \Cinv \ps h^{-1} \Vhnorm{\trialfunch[1]} \leq \Cinv^2 \ps^2 h^{-2} \Hhnorm{\trialfunch[1]},
 \end{equation}
 see, e.g., \cite{ErnG21a}*{Lemma 12.1}. 
 For a function $\trialfunc[1] \in \H2$ we moreover use, that 
 \begin{equation} \label{eq::RitzIntDiff}
	\bigVhnorm{(\NodalInt- \Rproj) \trialfunc[1]} \leq \CRitz h \norm{\trialfunc[1]}_{\H2} ,
 \end{equation}
 see, e.g., \cite{BreS08}*{Theorem 4.4.20 and 8.5.3}.
 Lastly, for all $\trialfunc[2] \in \H2$ it holds that $\trialfunch[2] = \NodalInt \trialfunc[2]$ satisfies 
 \begin{equation}
	\label{eq::InterpolationH1stable}
	\Hhnorm{\dLapl^{\half} \trialfunch[2]} = \Vhnorm{\NodalInt \trialfunc[2]} \leq \CInt \norm{\trialfunc[2]}_{\H2}	,
 \end{equation}
 see, e.g., \cite{ErnG21a}*{Section 11.5.1}. Note here that we could also require $\trialfunc[2] \in W^{1,s}(\fulldomain)$ with $s>m$, since we only need an embedding into $ C^{0}(\fulldomain)$ to make $\NodalInt \trialfunc[2]$ well defined, see, e.g. \cite{BreS08}*{Theorem 4.4.4}.

\begin{lemma}%
	\label{lem::dLaplonCNsol}
	Let $ \qExakt[0],\pExakt[0] \in  D(\mLapl) $  
	and $ f \in L^\infty\bigl([0,T]; D(\mLapl)\bigr)$ 
	. With $\qCN{0}=\qComponent_h^0$ and $\pCN{0}=\pComponent_h^0$ from \eqref{eq:projected-data}, the \CN approximations \eqref{eq::CN_on_VH_operators} satisfy
	\begin{equation}
		\label{dLaplonCNsol}
		\Bigl(  \bigVhnorm[]{\dLapl^{\half} \qCN{n-1}}^2 + \bigHhnorm[]{\dLapl^{\half} \pCN{n-1}}^2 \Bigr)^{\half} + \BigHhnorm[]{\ts[2] \dLapl^{1/2} \big( \fmh{n} \big)} \leq 	\Cdata,
	\end{equation}
	with $ \Cdata = \Cdata(\qExakt[0],\pExakt[0],\ps,f, t_n) $ independent of $h$.
\end{lemma}

\begin{proof}
	By the discrete variation-of-constants formula, we have
	\begin{equation*}
		\pmat{\qCN{n-1} \\ \pCN{n-1}} = \RCN^{n-1} \pmat{\qCN{0}\\ \pCN{0}} + \ts \sum_{j=1}^{n-1} \RCN^{n-j} \Rminus^{-1} \pmat{0 \\ \fmh{j}}.
	\end{equation*}
	It is well known that
	\begin{equation*}
		\bignorm{\RCN}_{\dEspace[]\leftarrow\dEspace[]} = 1 
		\quad \text{and} \quad 
		\bignorm{\Rminus}_{\dEspace[]\leftarrow\dEspace[]} \leq 1
	\end{equation*}
	so that we obtain
	\begin{equation*}
		\BigdEspacenorm{\pmat{\dLapl^{\half} \qCN{n-1} \\ \dLapl^{\half} \pCN{n-1}}} \leq  \BigdEspacenorm{\pmat{\dLapl^{\half} \qCN{0}\\ \dLapl^{\half} \pCN{0}}} + \ts \sum_{j=0}^{n-1} \BigHhnorm{\dLapl^{\half} f_h^j}.
	\end{equation*}
	It remains to bound $\Vhnorm{\dLapl^{\half} \qCN{0}}$, $\Hhnorm{\dLapl^{\half} \pCN{0}}$, and $\Hhnorm{\dLapl^{\half} f_h^j}$ in terms of the initial values and the inhomogeneity. Recall that the Ritz projection $\Rproj$ satisfies
	\begin{equation}
		\label{eq::RitzProjPropCont}
		\bigHhnorm{\dLapl \Rproj \phi} =  \bignorm{ \mLapl \phi}_{\L2},
	\end{equation} 
	for all $\phi \in D(L)$, since
		\begin{equation*}
		\ipML{\fulldomain}{\dLapl \Rproj \phi , \psi_h} = \ba{\fulldomain} \ip{\Rproj \phi,  \psi_h} =  \ba{\fulldomain} \ip{ \phi,  \psi_h} =\ipD{\fulldomain}{ \mLapl \phi, \psi_h} 
		, \quad \text{for all } \psi_h \in \Hh. 
	\end{equation*}
	Hence, we have from \eqref{eq:inverse-ineq} that
	\begin{equation*}
		\begin{aligned}
			\bigVhnorm{\dLapl^{\half}  \qCN{0}}
			 = \bigHhnorm{\dLapl \NodalInt \qExakt[0]}                                                                      %
			& \leq \bigHhnorm{\dLapl (\NodalInt - \Rproj) \qExakt[0]} + \bigHhnorm{\dLapl \Rproj \qExakt[0]}                 \\
			& \leq \Cinv h^{-1} \ps \bigVhnorm{(\NodalInt- \Rproj) \qExakt[0]} + \bignorm{ \mLapl \qExakt[0]}_{\L2}			 \\
			& \leq (\Cinv \CRitz +1) \ps \bignorm{\qExakt[0]}_{\H2},
		\end{aligned}
	\end{equation*}
	since $ \qExakt[0] \in D(\mLapl) \subset \H2 $ by assumption and \eqref{eq::RitzIntDiff}.
The remaining bounds in the $\Hh$-norm follow directly from \eqref{eq::InterpolationH1stable} and using that $ f \in L^\infty([0,T]; D(\mLapl))$. Setting all these bounds together yields an upper bound for \eqref{dLaplonCNsol} given by
\begin{equation}
	\label{eq::def_Cdata}
	\Cdata \! = \! \bigl( \Cinv \CRitz +1 \bigr)\ps\norm{\qExakt[0]}_{\H2} +  \CInt \ps \norm{\pExakt[0]}_{\H2} + t_n  \CInt \ps \! \max\limits_{j=0,\dots,n} \! \norm{f(t_j)}_{\H2}\;  .
\end{equation}
\end{proof}

Based on this stability result, we can bound the prediction error $\bigdEspacenorm[]{\solCNmod{n} - \solLFmod{n}}^2 $ which is \textit{global in space} but \textit{local in time} (we only consider one time step of the \CN and the \LF method starting from the \ds approximation $\solDS{n-1}$).

\begin{lemma}
	\label{PredictionErrorEstimate}
    Let the assumptions of \Cref{lem::dLaplonCNsol} and in addition
    the CFL condition \cref{eq::CFLconditionPrediction} be satisfied.
	Then it holds
	\begin{equation*}
		\bigdEspacenorm[]{\solCNmod{n} - \solLFmod{n}}  \leq \MaxCFLellc\left( \bigdEspacenorm{\solDS{n-1} - \solCN{n-1}} + \ts \Cdata \right),
	\end{equation*}
	with a constant $ \MaxCFLellc $ only depending on $ \ovp $.
\end{lemma}

\begin{proof}
	We rewrite the difference we are interested in with \Cref{def::operatorNotation} by
	\begin{equation}
		\label{eq::predictionDiff}
		\begin{aligned}
			\solCNmod{n} - \solLFmod{n}
			&= \big( \RCN - \RLF \big) \solDS{n-1} + \ts \big( \Rminus^{-1} - \RminusLF^{-1} \big) \gmh{n} \\
			&= \big( \RCN - \RLF \big) \big(\solDS{n-1} -\solCN{n-1} \big)
			+ \big( \RCN - \RLF \big)\solCN{n-1}
			+ \ts \big( \Rminus^{-1} - \RminusLF^{-1} \big) \gmh{n}
		\end{aligned}
	\end{equation}
	 with $\gmh{n} = \pmat{0 \\ \fmh{n}}$. A simple calculation shows
	 \begin{align*}
		\Rminus^{-1} - \RminusLF^{-1}
		& = \tss[4]\pmat{\dsysMat^{-1} \dLapl & \qquad \\ \qquad  & \dsysMat^{-1} \dLapl} 
		\pmat{\tss[4] \dLapl &  - \ts[2] \id \\[0.2em] \ts[2]\dLapl &   - \id} ,
	\end{align*} 
	and
	\begin{align*}
		\RCN - \RLF
		= \tss[4]\pmat{\dsysMat^{-1} \dLapl^2&  \qquad \\ \qquad  & \dsysMat^{-1} \dLapl^2 }
		\pmat{\tss[2]\id & \frac{\ts^3}{4} \id \\[0.2em] -\ts \id & \tss[2] \id} ,
	\end{align*}
	where we denote $\dsysMat = \id + \tss[4] \dLapl$. 
	Since for $\alpha\in [0,2]$, there holds
	\begin{equation*}
		0 \leq \frac{\xi^{\alpha / 2}}{1 + \xi} \leq 1, \qquad \text{for all } \xi \geq 0,
	\end{equation*}
	it follows that 
	\begin{equation*}
		\tss \bignorm{\dsysMat^{-1} \dLapl}_{\Hh \leftarrow \Hh} \leq 4 \quad \text{ and } \quad \tss \bignorm{\dsysMat^{-1} \dLapl}_{\Vh \leftarrow \Vh} \leq 4.
	\end{equation*}
	From this and \eqref{eq::CFLconditionPrediction} we can derive for $\solutionVector = \pmatT{\qComponent , \pComponent}$
	\begin{equation}
		\label{eq::estimateRmRhat-a}
		\begin{aligned}
		\bigdEspacenorm[]{\big( \RCN - \RLF \big) \solutionVector }^2
		\! =&   \BigVhnorm[]{\tss[4] \dsysMat^{-1} \dLapl^2 
		\bigl( \tss[2] \qComponent + \frac{\ts^3}{4}\pComponent \bigr) }^2 \! + \BigHhnorm[]{\tss[4] \dsysMat^{-1} \dLapl^2 
		\bigl( - \ts \qComponent+ \tss[2] \pComponent \bigr)  }^2  \\
		\leq&   2 \Bigl( \bigVhnorm[]{\tss[2] \dLapl \qComponent }^2 \! + \bigVhnorm[]{\ts \dLapl^{\half} \qComponent }^2 \! + \bigHhnorm[]{\ts \dLapl^{\half} \pComponent}^2 \! + \bigHhnorm[]{\tss[2] \dLapl \pComponent}^2\Bigr)  \\
		\leq&  2 (\ell^2 + 1) \Bigl( \bigVhnorm[]{\ts \dLapl^{\half}\qComponent}^2 + \bigHhnorm[]{\ts \dLapl^{\half}\pComponent}^2 \Bigr)\\
		=&  \MaxCFLellc^2 \bigdEspacenorm[]{\solutionVector }^2
		\end{aligned}
	\end{equation}
	with $\MaxCFLellc^2= 8(\ell^4 +\ell^2)$ being a polynomial of degree $4$ only depending on $\ovp$. Alternatively, we can also keep a prefactor $\tau$ by using 
	\begin{equation}
		\label{eq::estimateRmRhat-b}
		\bigdEspacenorm[]{\big( \RCN - \RLF \big)\solutionVector}^2 \leq 2(\ell^2 + 1) \ts \Bigl( \bigVhnorm[]{\dLapl^{\half} \qComponent }^2 + \bigHhnorm[]{\dLapl^{\half} \pComponent}^2  \Bigr)^{\half} .
	\end{equation}
	Analogously, we derive for $g = \pmatT{0 , f}$
	\begin{equation}
		\label{eq::estimateRmRhat-c}
		\begin{aligned}
		\bigdEspacenorm[]{ \big( \Rminus^{-1} - \RminusLF^{-1} \big) g}^2 &\leq  \bigHhnorm[]{\frac{\ts^3}{8} \dsysMat^{-1} \dLapl^{3/2} f}^2 + \bigHhnorm[]{\tss[4]\dsysMat^{-1} \dLapl f}^2  \leq 2 \bigHhnorm[]{\ts[2] \dLapl^{\half} f}^2.
		\end{aligned}
	\end{equation}

	From \eqref{eq::predictionDiff} we deduce
	\begin{equation}
		\label{eq::predicitionDiffNorms}
		\begin{aligned}
		\bigdEspacenorm[]{\solCNmod{n} - \solLFmod{n}} 
		\leq& \;  \bigdEspacenorm[]{\big( \RCN - \RLF \big) \big(\solDS{n-1} -\solCN{n-1} \big)}
		+ \bigdEspacenorm[]{\big( \RCN - \RLF \big)\solCN{n-1}}  \\
		& \; + \ts \bigdEspacenorm[]{ \big( \Rminus^{-1} - \RminusLF^{-1} \big) \gmh{n}} . 
	\end{aligned}
	\end{equation}
	For the first term we use \eqref{eq::estimateRmRhat-a} to get 
	\begin{equation*}
		\bigdEspacenorm[]{\big( \RCN - \RLF \big) \bigl( \solDS{n-1}- \solCN{n-1} \bigr) } \leq \MaxCFLellc \bigdEspacenorm[]{\solDS{n-1}- \solCN{n-1}},
	\end{equation*}
	while we take \eqref{eq::estimateRmRhat-b} for the second term to obtain
	\begin{equation*}
		\bigdEspacenorm[]{\big( \RCN - \RLF \big)\solCN{n-1}} \leq \sqrt{2 (\ell^2 +1)}\ts \Bigl( \bigVhnorm[]{\dLapl^{\half} \qCN{n-1} }^2 + \bigHhnorm[]{\dLapl^{\half} \pCN{n-1}}^2  \Bigr)^{\half},
	\end{equation*}
	and \eqref{eq::estimateRmRhat-c} to bound
	\begin{equation*}
		\bigdEspacenorm[]{ \big( \Rminus^{-1} - \RminusLF^{-1} \big) \gmh{n}} \leq \sqrt{2} \bigHhnorm[]{\ts[2] \dLapl^{\half} \fmh{n}}.
	\end{equation*} 
	To simplify the constants we just take 
		\begin{equation*}
			\max \{ \MaxCFLellc, \sqrt{2(\ell^2+1)}, \sqrt{2}\} = \MaxCFLellc, \qquad \forall \ell \in \N.
		\end{equation*}
	Putting these estimates together in \eqref{eq::predicitionDiffNorms} and using \Cref{lem::dLaplonCNsol} proves the claim.
\end{proof}

\section{Local deviation from global \CN approximation}
\label{sec::LocalError} 

Next, we investigate the local \CN approximations on all subdomains. After the for loop in Algorithm~\ref{algorithm::DSAlgWave} we end up with subdomain approximations $ \solDSloc{i}{n} $ on each of the overlapping subdomains $ \SDov{i} $, which we now compare to the global in space \CN approximation $\solCNmod{n}$ given in \eqref{FullDiscrCN-ML} on $ \fulldomain $. 
We refer to \Cref{table::notation:approx} for an overview on the notation.

\subsection{Variational problem for local error on subdomains}

Recall that the subdomain approximation  $ \solDSloc{i}{n} $ and the \CN approximation $\solCNmod{n}$  both start from the same approximation $ \solDS{n-1} $.
\begin{lemma}
	\label{lem::LocalDifferenceOnSubdomains}
	The difference $ \qDiffDSCN = \qDSloc{i}{n} - \restr{\qCNmod{n}}{\SDov{i}}$ satisfies
	\begin{subequations}
		\begin{equation}
			\label{eq::EllipticProblemLocalDifference}
			\begin{aligned}
				\ipML{\SDov{i}}{\qDiffDSCN, \varphi} + \tss[4] \baSDov{i}\ip{ \qDiffDSCN,  \varphi} & = 0,  &&\text{for all } \varphi \in \Vhn[\SDov{i}], \\
				\qDiffDSCN                                                                          &  %
				= \qLFmod{n} - \qCNmod{n}, \qquad && \text{on } \dSDov{i}.
			\end{aligned}
		\end{equation}
		Moreover, for $ \pDiffDSCN = \pDSloc{i}{n} - \restr{\pCNmod{n}}{\SDov{i}} $ it holds, that
		\begin{equation}
			\label{eq::EllipticProblemLocalDifference_p}
			\qDiffDSCN = \ts[2] \pDiffDSCN \quad \text{ in } \Hh[\SDov{i}].
		\end{equation}
	\end{subequations}
\end{lemma}
    The problem \eqref{eq::EllipticProblemLocalDifference} is well posed, cf.\ \cite{ErnG21b}*{Lemma 33.1}.

\begin{proof}
	The proof is based on comparing two weak forms, one on the overlapping subdomains $ \SDov{i} $ and one on the full domain $ \fulldomain $. 
	
	On the one hand, $  \solDSloc{i}{n}  \in \left( \Vh[\SDov{i}] \cap \Hn{1} \right) \times\Hh[\SDov{i}] $ is the solution of
	\begin{alignat*}{3}
		\nonumber \ipML{\SDov{i}}{\qDSloc{i}{n}, \varphi_i} -\ts[2]\ipML{\SDov{i}}{\pDSloc{i}{n}, \varphi_i} & = \ipML{\SDov{i}}{\qDS{n-1}, \varphi_i}      + \ts[2]\ipML{\SDov{i}}{\pDS{n-1}, \varphi_i}, \\
		\ipML{\SDov{i}}{\pDSloc{i}{n}, \psi_i} + \ts[2]\ba{\SDov{i}}\ip{\qDSloc{i}{n}, \psi_i}               & = \ipML{\SDov{i}}{\pDS{n-1}, \psi_i}         - \ts[2] \ba{\SDov{i}}\ip{\qDS{n-1}, \nonumber \psi_i}  + \ts \ipML{\SDov{i}}{\fmh{n}, \psi_i},    \\
		\nonumber \qDSloc{i}{n}                                                                              & = \qLFmod{n}, \qquad \text{on } \SDovint{i}. 
	\end{alignat*}
	for all $\varphi_i,\psi_i \in \Vhn[\SDov{i}]$.
	On the other hand, $  \solCNmod{n} \in \Vhn[\fulldomain]\times \Hhn[\fulldomain] $ satisfies
	\begin{alignat*}{3}
		\nonumber \ipML{\fulldomain}{\qCNmod{n}, \varphi} -\ts[2]\ipML{\fulldomain}{\pCNmod{n}, \varphi} & = \ipML{\fulldomain}{\qDS{n-1}, \varphi}  + \ts[2]\ipML{\fulldomain}{\pDS{n-1}, \varphi}, \\
		\ipML{\fulldomain}{\pCNmod{n}, \psi} + \ts[2]\ba{\fulldomain}\ip{\qCNmod{n}, \psi}               & = \ipML{\fulldomain}{\pDS{n-1}, \psi}    - \ts[2] \ba{\fulldomain}\ip{\qDS{n-1}, \nonumber \psi} + \ts \ipML{\fulldomain}{\fmh{n}, \psi},
	\end{alignat*}
	for all $ \varphi, \psi \in \Vhn[\fulldomain]$.

	In order to compare these two approximations, we extend the test functions $\varphi_i,\psi_i$ by zero onto $ \fulldomain $. By this extension we get valid test functions for the weak formulation of the global \CN on $ \fulldomain $, which vanish outside of $ \SDov{i} $. 
	The difference
	\begin{equation*}
		\solDSloc{i}{n} - \restr{\solCNmod{n}}{\SDov{i}} = \pmat{\qDSloc{i}{n}- \restr{\qCNmod{n}}{\SDov{i}} \\ \pDSloc{i}{n}- \restr{\pCNmod{n}}{\SDov{i}}} = \pmat{\qDiffDSCN \\[1ex] \pDiffDSCN}\in \Vh[\SDov{i}] \times \Hh[\SDov{i}]
	\end{equation*}
	satisfies
	\begin{subequations}
		\begin{align}
			\label{eq::DiffWeakForm-a}
				\ipML{\SDov{i}}{\qDiffDSCN, \varphi_i} -\ts[2]\ipML{\SDov{i}}{\pDiffDSCN, \varphi_i} & = 0,   & \varphi_i \in \Vhn[\SDov{i}], \\
			\label{eq::DiffWeakForm-b}
				\ipML{\SDov{i}}{\pDiffDSCN, \psi_i} + \ts[2]\ba{\SDov{i}}\ip{\qDiffDSCN, \psi_i}     & = 0,   & \psi_i \in \Vhn[\SDov{i}],
		\end{align}
		and
		\begin{equation}
			\label{eq::DiffWeakForm-c} 
			\begin{aligned}
				\restr{(\qDSloc{i}{n}-\qCNmod{n})}{\SDovint{i}}                 & = \restr{\qLFmod{n}}{\SDovint{i}} - \restr{\qCNmod{n}}{\SDovint{i}}, \\
				\restr{(\qDSloc{i}{n}-\qCNmod{n})}{\dfulldomain \cap \dSDov{i}} & = 0.
			\end{aligned}
		\end{equation}
	\end{subequations}
	\cref{eq::DiffWeakForm-a} is equivalent to \cref{eq::EllipticProblemLocalDifference_p}. Inserting this relation into \cref{eq::DiffWeakForm-b} and using the boundary conditions \eqref{eq::DiffWeakForm-c} shows that $ \qDiffDSCN  $ solves  \cref{eq::EllipticProblemLocalDifference}.
\end{proof}

\subsection{Exponential decay of errors in boundary data}
The goal of this subsection is to prove a discrete exponential decay result for the solution of \eqref{eq::DiffWeakForm-a}. For that we additionally need to introduce weighted bilinear forms. 
\begin{definition}{}
    \label{def::bilinearforms}
    Let $\lambda > 0$ be arbitrarily chosen and $\trialfunch[1],\trialfunch[2] \in \Vh$. 
    For a positive and continuous weight function  $ \weight: \fulldomain \to \R_+ $, we introduce the  weighted bilinear forms
    \begin{subequations}
        \begin{align}
            \baw{\D}{\weight}\ip{\trialfunch[1],\trialfunch[2]} 
            &= \int_{\D} \weight(x) \ps^2 \grad \trialfunch[1](x) \cdot \grad \trialfunch[2](x) \dintx , \label{def::baw}\\
            \bbw{\D}{\weight} \ip{\trialfunch[1],\trialfunch[2]} &= \frac1{\lambda^2}\ipML{\D}{\weight \trialfunch[1],\trialfunch[2]} + \baw{\D}{\weight}\ip{\trialfunch[1],\trialfunch[2]}\; . \label{def::bbw}
        \end{align}
        For  $\weight \equiv 1$ we abbreviate this via
        \begin{equation}
            \bb{\D}\ip{\trialfunch[1], \trialfunch[2]} = \bbw{\D}{1} \ip{\trialfunch[1],\trialfunch[2]} = \frac1{\lambda^2}\ipML{\D}{\trialfunch[1], \trialfunch[2]} + \ba{\D}\ip{\trialfunch[1], \trialfunch[2]} \label{def::bb}.
        \end{equation}
    \end{subequations}

\end{definition}

Next, we show a discrete decay estimate, which serves as key ingredient of our error analysis. This estimate is strongly motivated by \cite{BluLR92}*{Lemma 1}, which, however, does not apply directly to our setting for wave equations. For instance, we have to include mass lumping.
\begin{theorem} 
	\label{ExponentialDecay}
	Let 
	$ \lambda > 0 $, and $ g \in \Vh[\SDov{i}] $. 
	If $ z \in \Vh[\SDov{i}]$ solves
	\begin{equation}
		\label{eq::EllipticProblem}
		\begin{aligned}
			\ipML{\SDov{i}}{z, \varphi} + \lambda^2 \baSDov{i}\ip{ z,  \varphi} & = 0, \quad &&\forall \varphi \in \Vhn[\SDov{i}], \\
			z                                                                   & = g, \quad &&\text{on } \dSDov{i},
		\end{aligned}
	\end{equation}
	then it holds
	\begin{equation*}
		\frac{1}{\lambda^2}\ipML{\SD{i}}{z,z} + \baSD{i}\ip{z,z} \leq \beta \exp\Big(- \frac{\gamma \ov}{\max \{  \lambda, h \}}\Big) \left(\frac{1}{\lambda^2}\ipML{\SDov{i}}{g,g} +  \baSDov{i}\ip{g,g} \right)
	\end{equation*}
	with $ \gamma \in (0,1]$ and $ \beta > 0 $, which are independent of $ \lambda $, $ h, $ and the size of the subdomain $ \SD{i} $.
\end{theorem}

\begin{proof}
	We choose some $ \gamma \in (0,1] $, which will be fixed later and define a weight function
	\begin{equation}
		\label{eq::DefWeightFunction}
		\weight(x)  \coloneqq \exp\left(\frac{\gamma x}{\max \{ \lambda , h\}}\right), \quad x \in \SDov{i}.
	\end{equation}
	The distance of a point $ x \in \SDov{i} $ to the prediction interface is denoted by $ \dist(x, \SDovint{i})$. Since the distance function is not contained in $\Vh[\SDov{i}]$, 
	we define its nodal interpolation  $\NodalInt$ on the grid points in $\SDov{i}$ as
	\begin{equation}
		\label{eq::DefInterpolatedDistance}
		\Ipd \coloneqq \NodalInt \big( x \mapsto \dist(x, \SDovint{i}) \big) \in \Vh[\SDov{i}].
	\end{equation}
	For the sake of presentation, we suppress the space dependencies in the rest of the proof, i.e., we write  $ \weight(\Ipd) = \weight(\Ipd(x)) $.
	Let us first gather some  properties of $ \weight $ and $\Ipd$
	\begin{subequations}
		\begin{align}
			\weight(\xi) \geq 1 &\geq \weight(-\xi) \geq 0,  &&\text{ for }  \xi \geq 0, \label{eq::weightProps-a}\\
			\weight(\xi+\nu)& = \weight(\xi) \weight(\nu),  &&\text{ for } \xi,\nu \in \R, \label{eq::weightProps-b} \\
			\Ipd(x)- \ov &\geq 0 , &&\text{ for }  x  \in \SDc{i}.\label{eq::weightProps-c}
		\end{align}
	\end{subequations}
	By definitions  \eqref{def::bb} for $\D = \SD{i}$, \eqref{eq::DefWeightFunction}, \eqref{eq::weightProps-a}, and \eqref{eq::weightProps-c},  we have
	\begin{align}
		\label{eq::ExtensionWithWeights}
		\nonumber\bb{\SD{i}}\ip{z,z}
		          & = \frac{1}{\lambda^2} \ipML{\SD{i}}{z,z} + \baSD{i}\ip{z,z}                                      \\
		\nonumber & \leq \frac{1}{\lambda^2} \ipML{\SD{i}}{\weight(\Ipd-\ov)z,z} + \baw{\SD{i}}{\weight(\Ipd-\ov)}\ip{z,z} \\
		\nonumber & = \weight(-\ov) \bbw{\SD{i}}{\weight(\Ipd)}\ip{z,z}                                                 \\
		          & \leq \weight(-\ov) \bbw{\SDov{i}}{\weight(\Ipd)}\ip{z,z}.
	\end{align}
	Thus, it remains to bound $ \bbw{\SDov{i}}{\weight(\Ipd)}\ip{z,z} $ in terms of the boundary data $ g $. 
	To shorten the notation, we write  
	\begin{equation}
		\label{def::wz}
		\wz = \weight(\Ipd) z .
	\end{equation}
	Since $ \Iwz - g  \in \Vhn[\SDov{i}]$ is a valid test function for \eqref{eq::EllipticProblem}, we conclude that
	\begin{equation}
		\label{eq::Trick2}
			\frac{1}{\lambda^2}\ipML{\SDov{i}}{z, \Iwz} +  \baSDov{i}\ip{ z,  \Iwz  } = \frac{1}{\lambda^2}\ipML{\SDov{i}}{z, g } +  \baSDov{i}\ip{ z,  g } 
			=  \bbSDov{i}\ip{z,g} .
	\end{equation}
	Using this together with the product rule
	\begin{equation*}
		\weight(\Ipd) \grad z = \grad \wz - z \grad \weight(\Ipd) \, .
	\end{equation*}
	yields
		{\allowdisplaybreaks
			\begin{align}
				\nonumber \bbw{\SDov{i}}{\weight(\Ipd) }\ip{z,z}
				\nonumber & =                                                                       \frac{1}{\lambda^2} \ipML{\SDov{i}}{\wz,z} + \ba{\SDov{i}}\ip{\wz,z} - \int_{\SDov{i}} z \ps^2 \grad \weight(\Ipd) \cdot \grad z \dintx      \\\allowdisplaybreaks
				 &=  
				\frac{1}{\lambda^2} \ipML{\SDov{i}}{\wz - \Iwz,z}
				+  \ba{\SDov{i}}\ip{\wz- \Iwz,z}  
				+  \bbSDov{i}\ip{g,z}
				\label{eq::weightedCombinedNorm}
				\\
				& \quad
				- \int_{\SDov{i}} z \ps^2 \grad \weight(\Ipd) \cdot \grad z \dintx             
				\nonumber    
			\end{align}
		}
		The first term vanishes because of \eqref{eq::MLprojectIsInterpolation}. The third term can be estimated with \CSI by
		\begin{equation*}
			\abs{	\bb{\SDov{i}}\ip{g,z}} 
            \leq \bb{\SDov{i}}\ip{g,g}^{\half} \bb{\SDov{i}}\ip{z,z}^{\half} 
            \leq \bb{\SDov{i}}\ip{g,g}^{\half} \bbw{\SDov{i}}{ \weight(\Ipd) }\ip{z,z}^{\half}, 
		\end{equation*}
		where we also used that $ \weight(\Ipd) \geq 1 $ on $ \SDov{i} $.
		The bounds on the remaining two terms are given in Lemmas~\ref{Lem::BoundDecayTerm2} and \ref{Lem::BoundDecayTerm3} which we postponed to \Cref{sec::Appendix} 
		because the proofs are rather technical. 
		Plugging these bounds into $\eqref{eq::weightedCombinedNorm}$ yields
		\begin{equation}
			\label{eq::absorbtionrdy}
			\bbw{\SDov{i}}{\weight(\Ipd) }\ip{z,z}^{\half} \leq \gamma (C_1 + C_2) \bbw{\SDov{i}}{\weight(\Ipd) }\ip{z,z}^{\half} + \bb{\SDov{i}}\ip{g,g}^{\half}
		\end{equation} 
		with constants $C_1, C_2$ independent of $\gamma, h $ and $\lambda$.
		Thus, by possibly reducing $ \gamma \in (0,1] $, such that
		\begin{equation*}
		\beta^{-\half} \coloneqq	1 - \gamma \Big( C_1 + C_2 \Big) > 0 .
		\end{equation*}
		This allows us to rearrange \eqref{eq::absorbtionrdy}, which gives
		\begin{equation}
			\label{eq::BoundByBoundaryData}
			\begin{aligned}
				 \bbw{\SDov{i}}{\weight(\Ipd) }\ip{z,z}^{\half} & \leq \beta^{\half} \bb{\SDov{i}}\ip{g,g}^{\half}  
			\end{aligned}
		\end{equation}
		By setting together \eqref{eq::ExtensionWithWeights} and \eqref{eq::BoundByBoundaryData} we finally obtain
		\begin{equation*}
			\bb{\SD{i}}\ip{z,z}  \leq \beta \weight(- \ov)  \bb{\SDov{i}}\ip{g,g},
		\end{equation*}
		which concludes the proof. \qedhere

\end{proof}

We now set together \Cref{lem::LocalDifferenceOnSubdomains} and \Cref{ExponentialDecay} to apply the decay result to our setting.

\begin{lemma}
	\label{cor::BoundLocalDiff}
	Let the conditions from \Cref{lem::LocalDifferenceOnSubdomains} and \Cref{ExponentialDecay} be satisfied. Then,
	\begin{equation*}
		\bigdEspacenorm[\SD{i}]{\solDSloc{i}{n} - \solCNmod{n}}^2
		\leq \beta \exp\Big(-\frac{\gamma \ov}{\max\{\ts[2], h\}}\Big) \bigdEspacenorm[\SDov{i}]{\solLFmod{n}-\solCNmod{n}}^2 \; .
	\end{equation*}
\end{lemma}

\begin{proof}
	As in \Cref{lem::LocalDifferenceOnSubdomains}, we write 
	\begin{equation*}
		z_i^n = \solDSloc{i}{n} - \restr{\solCNmod{n}}{\SDov{i}} = \pmat{\qDSloc{i}{n}- \restr{\qCNmod{n}}{\SDov{i}} \\ \pDSloc{i}{n}- \restr{\pCNmod{n}}{\SDov{i}}} = \pmat{\qDiffDSCN \\[1ex] \pDiffDSCN}.
	\end{equation*}
	From \eqref{eq::EllipticProblemLocalDifference_p} it holds
	\begin{equation}
		\label{eq::dEspaceNorm-zin}
		\bigdEspacenorm[\SD{i}]{z_i^n}^2 = \bigVhnorm[\SD{i}]{\qDiffDSCN}^2 + \bigHhnorm[\SD{i}]{\pDiffDSCN}^2 = \bigVhnorm[\SD{i}]{\qDiffDSCN}^2 + \frac{4}{\tss}\bigHhnorm[\SD{i}]{\qDiffDSCN}^2 .
	\end{equation}
	The leapfrog steps \eqref{FullDiscrLF-ML-a} and \eqref{FullDiscrLF-ML-c} started from $\solDS{n-1}$ yield
\begin{align*}
    \qLFmod{n} 
            = \qDS{n-1} + \ts[2]\big( \pDS{n-1} + \pLFmod{n} \big)
\end{align*}
On the other hand, we have for the \CN method by \eqref{FullDiscrCN-ML-a} that 
\begin{equation*}
    \qCNmod{n} = \qDS{n-1} + \ts[2]\big( \pDS{n-1} + \pCNmod{n} \big).
\end{equation*}
It follows
\begin{equation*}
    \qLFmod{n}- \qCNmod{n} = \ts[2] \bigl( \pLFmod{n}- \pCNmod{n} \bigr),
\end{equation*}
and thus we have as in \eqref{eq::dEspaceNorm-zin}
\begin{equation*}
	\bigdEspacenorm[\SDov{i}]{\solLFmod{n}- \solCNmod{n}}^2  = 
    \bigVhnorm[\SDov{i}]{\qLFmod{n}- \qCNmod{n}}^2 
    + \frac{4}{\tss} \bigHhnorm[\SDov{i}]{\qLFmod{n}- \qCNmod{n}}^2   
    \; .
\end{equation*}
The stated bound then follows immediately by \Cref{ExponentialDecay}, where we set $g = \qLFmod{n}- \qCNmod{n} $ and $\lambda = \ts[2]$.
\end{proof}

\section{Averaging}
\label{sec::AvgProps}

In this section we investigate the stability of the averaging operator $\averaging$ defined in \eqref{eq::averaging:mean} and deduce some bound for the specific use of $\averaging$ within a single \ds step.
The statement of the following lemma relies on the notation summarized in \Cref{table::notation:approx}.

\begin{lemma}
	\label{lem::AvgProps}
	\begin{subequations}
		The averaging operator $ \averaging $ satisfies the following stability estimate in the mass lumped $ L^2 $-norm
		\begin{equation}
			\label{eq::avg_prop2}
			\bigHhnorm[\fulldomain]{\averaging \big( \big\{ \pComponent_i \big\}_{i=1,\dots,\numberSD}\big)}^2 \leq \cavg \sumSD \bigHhnorm[\SD{i}]{\pComponent_i}^2, \; \text{for all } \pComponent_i \in \Hh[\SDov{i}], \; i = 1,\dots \numberSD ,
		\end{equation}
		with a constant $ \cavg > 0 $ independent of $ h,\ps,\ell $.

		If the CFL condition \Cref{eq::CFLconditionPrediction} is satisfied, it further holds
		\begin{equation}
			\label{eq::avg_prop3}
			\begin{aligned}
				\bigdEspacenorm{\solDS{n}- \solCNmod{n}}^2
				 & = \bigdEspacenorm{\avOpDiff}^2                                                                 \\
				 & \leq \Cavg \sumSD \bigdEspacenorm[\SD{i}]{\solDSloc{i}{n}- \solCNmod{n}}^2 .
			\end{aligned}
		\end{equation}
		There, the constant $ \Cavg $ may depend on $ \ovp $, but it does not depend on $ h $ and $ \ts $. 
	\end{subequations}
\end{lemma}
Note that \eqref{eq::avg_prop3} relies heavily on the specific structure of $ \solDSloc{i}{n}- \restr{\solCNmod{n}}{\SDov{i}} $ from \Cref{lem::LocalDifferenceOnSubdomains}.

\begin{proof}
    With $\pComponent_{\fulldomain} = \averaging \big( \big\{ \pComponent_i \big\}_{i=1,\dots,\numberSD}\big) $ it holds, that
    \begin{equation*}
        \bigHhnorm[\fulldomain]{v_{\fulldomain}}^2
        = \sumSD \bigHhnorm[\SD{i}]{v_{\fulldomain}}^2
        = \sumSD \sum_{K\in \Th[\SD{i}]} \sum_{j=1}^{\dim+1} \frac{\abs{K}}{\dim+1} \pComponent_{\fulldomain}(\xcoord_{K,j})^2.
    \end{equation*}
    The stability bound \cref{eq::avg_prop2} follows directly by 
    \begin{equation*}
        \pComponent_{\fulldomain}(\xcoord_{K,j})^2 = \left(  \sum\limits_{k \in J_{K,j} } \frac{\pComponent_k(\xcoord_{K,j})}{\abs{J_{K,j}}} \right)^2 
		\leq \sum\limits_{k \in J_{K,j} } \frac{\pComponent_k(\xcoord_{K,j})^2}{\abs{J_{K,j}}}, 
    \end{equation*}
    with $J_{K,j} =  \big\{ k \in \{1,\dots, \numberSD \} : \xcoord_{K,j} \in \SDc{k}\big\}$ and the fact that the sizes of neighboring cells are uniformly equivalent w.r.t.\ $h$, see \cite{ErnG21a}*{Proposition 11.6}.

	To show \cref{eq::avg_prop3} we note that 
	\begin{equation*}
		\begin{aligned}
			\bigdEspacenorm{\avOpDiff}^2 
			= \sumSD \bigdEspacenorm[\SD{i}]{\avOpDiff}^2 . 
		\end{aligned}
	\end{equation*}
	Let $ K \in \Th$ be a cell of which all vertices $ \{\xcoord_{K,j}\}_{j=1}^{\dim+1} $ lie inside of one of these subdomains, say in $ \SD{k} $. For this cell the averaged solution coincides with the local solution, i.e.,
	\begin{equation}
		\label{eq::AvgTrivialAwayFromInterface}
		\restr{\avOpDiff}{K} = \restr{\Bigl(\solDSloc{k}{n} - \restr{\solCNmod{n}}{\SD{k}}\Bigr)}{\K} .
	\end{equation}
	This means we only have to focus on the remaining cells for which one of the vertices, say $\xcoord_{K,r}$, lies on an interface, i.e. in 
	\begin{equation*}
		\xcoord_{K,r} \in  \Avgint = \bigcup\limits_{i=1,\dots,\numberSD} \SDint{i}.
	\end{equation*}
	Furthermore, we define the set of all cells having at least one vertex on the interface by
	\begin{equation*}
		\label{def::CellsAroundInterfaces}
		\Th[\Avgint]  \! \coloneqq \! \big\{ K \in \Th \;|\; \xcoord_{K,j} \in \nodalPoints \setminus \big( \mathop{\bigcup}\limits_{i=1,\dots,\numberSD}  \nodalPoints[\SD{i}]\big) \text{ for at least one } j \in {1,\dots,\dim+1} \big\} .
	\end{equation*}
	By \Cref{lem::LocalDifferenceOnSubdomains} we have
	$
		\qDSloc{i}{n} - \qCNmod{n} = \ts[2] \left( \pDSloc{i}{n} - \pCNmod{n}\right)
	$
    in $\Hh[\SDov{i}]$, see \eqref{eq::EllipticProblemLocalDifference_p}, and thus 
	\begin{equation*}
		\left(\qDS{n} - \qCNmod{n}\right)(\xcoord_{K,j}) = \ts[2] \left( \pDS{n} - \pCNmod{n} \right) (\xcoord_{K,j}) \quad \text{for all }  \xcoord_{K,j}  \in \nodalPoints[] \; .
	\end{equation*}
    For $K \in \Th[\Avgint]$, we derive with the CFL condition \eqref{eq::CFLconditionPrediction} that
	\begin{align*}
		 \bigdEspacenorm[\K]{\solDS{n}- \solCNmod{n}}^2
		 & = \bigVhnorm[\K]{\qDS{n}- \qCNmod{n}}^2 + \frac4{\tss} \bigHhnorm[\K]{\qDS{n}- \qCNmod{n}}^2 \\
		 & \leq \bignorm{\dLapl}_{\Hh \leftarrow \Hh}
		\bigHhnorm[\K]{\qDS{n}-\qCNmod{n}}^2 + \frac4{\tss} \bigHhnorm[\K]{\qDS{n}- \qCNmod{n}}^2 \\
		& \leq(\ovp^2 + 1)\frac4{\tss} \bigHhnorm[\K]{\qDS{n}- \qCNmod{n}}^2 \\
		& =  (\ovp^2 + 1) \bigHhnorm[\K]{\pDS{n}- \pCNmod{n}}^2 \; .
	\end{align*}
	We continue by using \cref{eq::avg_prop2} which implies that
	\begin{align*}
		\sumTh[\Avgint] \bigHhnorm[\K]{\pDS{n}- \pCNmod{n}}^2
		&\leq \cavg \sumSD \sum\limits_{\substack{K\in \Th[\Avgint] \\ K \subset \SDc{i}}} \bigHhnorm[\K]{\pDSloc{i}{n}- \pCNmod{n}}^2 \\
		&\leq \cavg \sumSD \sum\limits_{\substack{K\in \Th[\Avgint] \\ K \subset \SDc{i}}} \bigdEspacenorm[\K]{\solDSloc{i}{n}- \solCNmod{n}}^2 .
	\end{align*}
	Putting these derivations together we obtain finally
	\begin{equation}
		\label{eq::AvgOnInterfaces}
		\sumTh[\Avgint] \bigdEspacenorm[\K]{\solDS{n}- \solCNmod{n}}^2 \leq \Cavg \sumSD \sum\limits_{\substack{K\in \Th[\Avgint] \\ K \subset \SDc{i}}} \bigdEspacenorm[\K]{\solDSloc{i}{n}- \solCNmod{n}}^2 \, ,
	\end{equation}
	with $ \Cavg =  \cavg \big(\ovp^2 + 1\big)   $. However, by \cref{eq::AvgTrivialAwayFromInterface} we also know that
	\begin{equation*}
		\sum_{K \in \Th \setminus \Th[\Avgint]} \bigdEspacenorm[\K]{\solDS{n}- \solCNmod{n}}^2 = \sumSD \sum\limits_{\substack{K\in \Th \setminus \Th[\Avgint] \\ K \subset \SDc{i}}} \bigdEspacenorm[\K]{\solDSloc{i}{n}- \solCNmod{n}}^2 \; .
	\end{equation*}
	Together with \cref{eq::AvgOnInterfaces} this yields \cref{eq::avg_prop3}.
\end{proof}

\section{Global error}
\label{sec::ErrorRecursion}
Recall the definitions of $\solDS{n}, \solDSloc{i}{n} ,\solCN{n}, \solCNmod{n}$ and $\solLFmod{n}$, which are summarized in \Cref{table::notation:approx}.
We now derive a bound for the difference $\errorDSCN{n}$ of the \ds approximation and the \CN approximation after $n$ \timesteps, i.e. 
\begin{equation}
	\label{def::errorDSCN}
	\errorDSCN{n} = \bigdEspacenorm[]{\solDS{n}-\solCN{n}} \; .
\end{equation}
For that we will use the results of the previous sections, so recall the constants $\Cdata$ from \Cref{lem::dLaplonCNsol}, $\MaxCFLellc$ from \Cref{PredictionErrorEstimate}, $\beta, \gamma$ from \Cref{ExponentialDecay}, and $\Cavg$ from \Cref{lem::AvgProps}.

To set everything together, we denote the maximum number of mutually overlapping subdomains by
\begin{equation}
	\label{const::Cglob}
	\card \{ i \; | \; x \in \SDov{i} \} \leq \Cglob, \qquad \text{for all }   x\in \fulldomain .
\end{equation}
Note that $\Cglob$ does not depend on the total number of subdomains $\numberSD$.
This property is known as \textit{bounded local overlap property}, see for instance \cite{Hol03}.

\begin{lemma}%
	\label{ErrorRecursion}
	Let $ \qExakt[0],\pExakt[0] \in  \H2 \cap \Hn1  $ and $ f \in L^{\infty}([0,T]; \H2 \cap \Hn1) $.
	Assume that the CFL condition \eqref{eq::CFLconditionPrediction} holds. 
	Then, $\errorDSCN{n}$ defined in \eqref{def::errorDSCN}
    satisfies 
    \begin{equation}
        \label{eq::errorRecursion}
        \errorDSCN{n} \leq   \ts \Cdata \Derror \sum_{j=0}^{n-1}  (1+\Derror)^j
    \end{equation}
    with 
    \begin{equation}
        \label{def::D-errorRecursion}
        \Derror^2 =   \Cavg \Cglob \beta \exp\Big(-\frac{\gamma \ov}{\max\{\ts[2], h\}}\Big) \MaxCFLellcs.
    \end{equation}
\end{lemma}

Note, that the overlap $ \ov $ scales like $ \ovp h $. Recall, that $\MaxCFLellcs$ is a fourth polynomial of order in $\ovp$. Hence, $\Derror$ becomes small fast when increasing $\ovp$.

\begin{proof}
	Firstly, we use the triangle inequality to obtain 
	\begin{equation*}
		\errorDSCN{n} = \bigdEspacenorm{\solDS{n} - \solCN{n}} \leq \bigdEspacenorm{\solDS{n} - \solCNmod{n}} + \bigdEspacenorm{\solCNmod{n} - \solCN{n}} = E^n_A + E^n_B.
	\end{equation*}
	Secondly, we use \eqref{eq::avg_prop3} and then exploit the exponential decay in \Cref{cor::BoundLocalDiff} and \eqref{const::Cglob}
	\begin{align*}
		(E^n_A)^2  =\bigdEspacenorm{\solDS{n} - \solCNmod{n}}^2
		&\leq \Cavg \sumSD \bigdEspacenorm[\SD{i}]{\solDSloc{i}{n} - \solCNmod{n}}^2 \\
		&\leq \Cavg  \beta \exp\Big(-\frac{\gamma \ov}{\max\{\ts[2], h\}}\Big) \sumSD \bigdEspacenorm[\SDov{i}]{\solLFmod{n}- \solCNmod{n}}^2 \\
		&\leq \Cglob \Cavg  \beta \exp\Big(-\frac{\gamma \ov}{\max\{\ts[2], h\}}\Big) \bigdEspacenorm{\solLFmod{n}- \solCNmod{n}}^2
	\end{align*}
	Next, \Cref{PredictionErrorEstimate} yields
	\begin{equation*}
		\bigdEspacenorm{\solLFmod{n} - \solCNmod{n}} \leq \MaxCFLellc\left( \errorDSCN{n-1} + \ts \Cdata \right) ,
	\end{equation*}
   so that we obtain
	\begin{align*}
		(E^n_A)^2 
		 & \leq \Cavg \Cglob \beta \exp\Big(-\frac{\gamma \ov}{\max\{\ts[2], h\}}\Big)
		 \MaxCFLellcs\left( \errorDSCN{n-1} + \ts \Cdata \right)^2 .
	\end{align*}
	For term $E_B^n$ we use that the \CN operator $ \RCN $ defined in \eqref{eq::def_R_operators} is unitary, which implies
	\begin{equation*}
		E^n_B = \bigdEspacenorm{\solCNmod{n} - \solCN{n}} = \bigdEspacenorm{\RCN\left(\solDS{n-1} - \solCN{n-1}\right)} =  \bigdEspacenorm{\solDS{n-1} - \solCN{n-1}} =  \errorDSCN{n-1} .
	\end{equation*}
    Combining the estimates of $E_A^n$ and $E_B^n$ and obtain
	\begin{align*}
		\errorDSCN{n} & \leq  \Derror \left( \errorDSCN{n-1} + \ts \Cdata   \right) + \errorDSCN{n-1} \\
                      & = (1+\Derror)\errorDSCN{n-1} + \tau \Derror \Cdata.
	\end{align*}
	Resolving this inequality and using that $
		\errorDSCN{0}
		= \bigdEspacenorm{\solDS{0} - \solCN{0}}
		= 0
	$ proves the claim.
\end{proof}

\subsection{Proof of \texorpdfstring{\Cref{ConvergenceDS}}{}}
\label{sec::ErrorAnalysis}
In this section, we complete the proof of the main result \Cref{thm::main} which relies on \Cref{ConvergenceDS}.  

\begin{proof}
	By assumption, we take the overlap $ \ov \sim \ovp h $ sufficiently large, such that there exists a  constant $ 0 < \sigma \leq  1 $ with
	\begin{equation}
		\label{eq::AssumptionOverlap}
		\Derror = \left(  \Cavg \Cglob \beta \exp\Big(-\frac{\gamma \ov}{\max\{\ts[2], h\}}\Big)\right)^{\half} \MaxCFLellc \leq \sigma \tss,
	\end{equation}
	where $\gamma$ was introduced in \Cref{ExponentialDecay}.
	From \eqref{eq::errorRecursion} we obtain
	\begin{equation}
		\label{eq::ErrorRecursion}
		\errorDSCN{n} \leq \ts^2  \Cdata\sigma \ts \sum_{j=0}^{n-1} \left( 1 + \sigma \tss \right)^j.
	\end{equation}
	This allows us to derive two error bounds.
	Since $\ts \leq 1$, there holds 
	\begin{align*}
		\sigma \ts \sum_{j=0}^{n-1} \left( 1 + \sigma \ts^2 \right)^j 
		&\leq \sigma \ts \sum_{j=0}^{n-1} \left( 1 + \sigma \ts \right)^j 
		\leq %
		\e^{\sigma t_n} - 1 \\
		\intertext{and}
		\sigma \ts \sum_{j=0}^{n-1} \left( 1 + \sigma \tss \right)^j 
		 &\leq  \frac1{\ts} \bigl( \e^{\sigma \ts t_n} -1 \bigr) 
		 \leq  \sigma t_n \e^{\sigma \ts t_n}.
	\end{align*}
	This proves the lemma.
\end{proof}

\subsection{Discussion on the CFL condition}
\label{sec::CFLdiscussion}

Next, let us discuss the CFL condition \eqref{eq::CFLconditionPrediction} required for \Cref{thm::main}. Obviously, \eqref{eq::CFLconditionPrediction} becomes weaker if $\ell$ is chosen larger.

The second condition is that the overlap $\ov \sim \ovp h$ has to be chosen sufficiently large such that \cref{eq::AssumptionOverlap} holds for a constant $\sigma \in (0,1)$. 
This condition ensures, that the damping of the prediction error is strong enough.  
The condition can be satisfied for fixed $ \ts, h $ by taking $ \ovp $ sufficiently large such that the condition on $\ov \sim h \ovp $ above is satisfied.

\section{Numerical experiments}
\label{sec::NumericalExp}

After proving the theoretical results in the previous sections, we confirm and illustrate these with some numerical experiments. All experiments were implemented within the FEniCSx framework \cites{Fenics, Fenics2}. The code corresponding to the experiments in this section is made publicly available at
\begin{equation*}
	\text{\url{https://doi.org/10.35097/zsd8jg6d7urs378x}}.
\end{equation*}

In all graphs, solid lines refer to errors against the analytical solution, dashed lines refer to errors against the \CN\ solution.

\subsection{One-dimensional experiments} \label{subsec:1dex}

Inspired by \cite{Lin22} we construct an analytical solution to the homogeneous wave equation on a one-dimensional domain $ (0,1) \subseteq \R  $.  We choose a function
\begin{equation} \label{eq:mu-xi-s}
	\mu_{\xi,s}(z) = \mathbbm{1}_{\{ \abs{z -\xi}<s\}} \sin\Bigl( \mfrac{ z -(\xi+s)}{2s}\pi \Bigr)^3 ,
\end{equation}
where $\mathbbm{1}_{S}$ denotes the indicator function corresponding to a set $S$.
Then, we define the initial conditions for $\eqref{eq::waveEqFirstOrder}$ as
\begin{equation} \label{eq:1d-initial}
	\qExakt[0](x) = \mu(x), \qquad \pExakt[0](x) = \restr{\dt \mu(x-t)}{t=0}
\end{equation}
with
\begin{equation*}
	\mu(z) = \mu_{\xi_1,s}(z) - \mu_{\xi_2,s}(z),
	\qquad 
	\xi_1 = 0.55, \quad \xi_2 = 0.45, \quad  s = 0.2. 
\end{equation*}

\begin{figure}
	\begin{tikzpicture}

	\definecolor{darkgray176}{RGB}{176,176,176}

	\begin{axis}[
		width = \textwidth,
		height = 0.35\textwidth,
			legend cell align={left},
			legend style={
					cells={align=left},
					fill opacity=0.8,
					draw opacity=1,
					text opacity=1,
					at={(0.03,0.97)},
					anchor=north west,
					draw=none
				},
			tick pos=both,
			x grid style={darkgray176},
			xlabel={\(\displaystyle h_{{\min}} \ell\)},
			xmin=1.53312890873946e-05, xmax=0.013,
			xtick style={color=black},
			xtick={0.002,0.005,0.008,0.01,0.012},
			y grid style={darkgray176},
			ymin=8e-05, ymax=0.0095,
			ytick style={color=black},
			ytick={0.002,0.004,0.006,0.008},
			ylabel={\(\displaystyle \tau_{\max} \)},
		]
		\addplot [KITblue, mark=*, mark size=2, mark options={solid,fill opacity=1}]
table {%
0.000306625781747893 0.0008111615833874
0.000613251563495787 0.0012916559028674
0.00091987734524368 0.0009769441187964
0.00122650312699157 0.0012916559028674
0.00153312890873947 0.0011233430689732
0.00183975469048736 0.0018740629685157
0.00214638047223525 0.0015556938394523
0.00245300625398315 0.0017070672584499
0.00275963203573104 0.0022573363431151
0.00306625781747893 0.0022573363431151
0.00337288359922683 0.0023640661938534
0.00367950938097472 0.0025947067981318
0.00398613516272261 0.0025947067981318
0.00429276094447051 0.0029832935560859
0.0045993867262184 0.0029832935560859
0.0049060125079663 0.0032743942370661
0.00521263828971419 0.0032743942370661
0.00551926407146208 0.0034293552812071
0.00582588985320998 0.0035945363048166
0.00613251563495787 0.0039432176656151
0.00643914141670576 0.0039432176656151
0.00674576719845366 0.0041322314049586
0.00705239298020155 0.0041322314049586
0.00735901876194944 0.0045330915684496
0.00766564454369734 0.0047528517110266
0.00797227032544523 0.0049751243781094
0.00827889610719312 0.0049751243781094
0.00858552188894102 0.0052137643378519
0.00889214767068891 0.0054644808743169
0.0091987734524368 0.0054644808743169
0.0095053992341847 0.005720823798627
0.00981202501593259 0.005720823798627
0.0101186507976805 0.0062814070351758
0.0104252765794284 0.006578947368421
0.0107319023611763 0.0068965517241379
0.0110385281429242 0.0068965517241379
0.0113451539246721 0.0072254335260115
0.01165177970642 0.0075642965204236
0.0119584054881678 0.0075642965204236
0.0122650312699157 0.0079239302694136
};
\addlegendentry{$\text{DS}_{\ell}$}
\addplot [gray, dotted, line width=1.2]
table {%
0.000306625781747893 0.000605582291607824
0.000613251563495787 0.000782400998244232
0.00091987734524368 0.000959219704880641
0.00122650312699157 0.00113603841151705
0.00153312890873947 0.00131285711815346
0.00183975469048736 0.00148967582478987
0.00214638047223525 0.00166649453142627
0.00245300625398315 0.00184331323806268
0.00275963203573104 0.00202013194469909
0.00306625781747893 0.0021969506513355
0.00337288359922683 0.00237376935797191
0.00367950938097472 0.00255058806460832
0.00398613516272261 0.00272740677124473
0.00429276094447051 0.00290422547788113
0.0045993867262184 0.00308104418451754
0.0049060125079663 0.00325786289115395
0.00521263828971419 0.00343468159779036
0.00551926407146208 0.00361150030442677
0.00582588985320998 0.00378831901106318
0.00613251563495787 0.00396513771769959
0.00643914141670576 0.00414195642433599
0.00674576719845366 0.0043187751309724
0.00705239298020155 0.00449559383760881
0.00735901876194944 0.00467241254424522
0.00766564454369734 0.00484923125088163
0.00797227032544523 0.00502604995751804
0.00827889610719312 0.00520286866415445
0.00858552188894102 0.00537968737079085
0.00889214767068891 0.00555650607742726
0.0091987734524368 0.00573332478406367
0.0095053992341847 0.00591014349070008
0.00981202501593259 0.00608696219733649
0.0101186507976805 0.0062637809039729
0.0104252765794284 0.00644059961060931
0.0107319023611763 0.00661741831724571
0.0110385281429242 0.00679423702388212
0.0113451539246721 0.00697105573051853
0.01165177970642 0.00714787443715494
0.0119584054881678 0.00732469314379135
0.0122650312699157 0.00750151185042776
};
\addlegendentry{$0.577h_{\min} \ell$}
	\end{axis}

\end{tikzpicture}
	\caption{Maximal stable step sizes $\tau_{\max}$ of the \ds\  (DS$_\ell$) depending on the overlap parameter $\ell$.}
	\label{Fig:1d_CFL}
\end{figure}

\begin{figure}
	
		\input{pictures/experiment_thin.tex}
		\caption{Errors of the \ds\ approximation (DS$_\ell$) for the 1d example of \Cref{subsec:1dex} for various numbers of  $\ell$ measured at the final time $ T=5.0 $ against the exact solution (solid) and against the \CN\ approximation (dashed). For comparison, we also show the error of the
			\CN\ (CN) and the \LF\ (lf) methods on the full domain.}
			\label{Fig:1D_ell_dependence}
	
\end{figure}

In order to avoid over approximation effects, we perturbed a uniform mesh randomly, which resulted in a non-equidistant mesh with minimum mesh width $ h_{\min} = 3\cdot 10^{-4} $ and maximum mesh width $ h_{\max} = 6.9\cdot 10^{-4} $ with $2000$ subintervals).
All \ds approximations use two subdomains and the index indicates the choice of the overlap parameter $ \ell $. 

We start by illustrating the dependence of the CFL condition \eqref{eq::CFLconditionPrediction} on the overlap parameter $\ell$
in \Cref{Fig:1d_CFL}. One can clearly see the linear dependence, confirming the theoretical result.

In \Cref{Fig:1D_ell_dependence} we show the $\Hn{1} \times \L{2}$ errors measured at the final time $ T=5.0 $ against the exact solution (which can be found in \cite{Lin22}). For comparison, we also show the error of \CN\ (CN) and the \LF\ (lf) methods on the full domain. Here, we observe that the error of DS$_\ell$ roughly coincides with the \CN\ error, if the time step satisfies the CFL condition.  This is confirmed in \Cref{Fig:1D_ell_dependence}, where we illustrate \Cref{ConvergenceDS} by plotting not only the error against the exact solution but also the error of the \ds\ method (DS$_8$) against the \CN\ solution on the full domain. 

\subsection{Two-dimensional experiments} \label{subsec:2dex}

Let $ \fulldomain = (0,1)^2$, $ \xi = 0.5 $,
$ s = 0.2 $, and
Using \eqref{eq:mu-xi-s}, we prescribe the exact solution of \eqref{eq::waveEqFirstOrder} as
\begin{equation*}
   \qComponent(x,y,t) = \qComponent_{\mathrm{1D}}(x,t)\mu_{\xi,s}(y) + \qComponent_{\mathrm{1D}}(y,t)\mu_{\xi,s}(x),
\end{equation*}  
where $u_{\mathrm{1D}}$ solves the homogeneous one-dimensional wave equation on $(0,1)$ with initial conditions \eqref{eq:1d-initial} with $\mu=\mu_{\xi,s}$. Inserting this solution into the wave operator defines the inhomogeneity $f$ in \eqref{eq::waveEqFirstOrder}.

\begin{figure}
	\begin{tikzpicture}

\definecolor{burlywood221204119}{RGB}{221,204,119}
\definecolor{darkgray176}{RGB}{176,176,176}
\definecolor{darkslateblue5134136}{RGB}{51,34,136}
\definecolor{forestgreen1711951}{RGB}{17,119,51}
\definecolor{gray}{RGB}{128,128,128}
\definecolor{indianred204102119}{RGB}{204,102,119}

\begin{axis}[
	width=0.96\textwidth,
	height=0.5\textwidth,
legend cell align={left},
legend style={
  fill opacity=0.8,
  draw opacity=1,
  text opacity=1,
  at={(0.97,0.03)},
  anchor=south east,
  draw=none
},
log basis x={10},
log basis y={10},
tick pos=both,
unbounded coords=jump,
x grid style={darkgray176},
xlabel={\(\displaystyle \tau\)},
xmin=0.000397164117362141, xmax=0.0629462705897084,
xmode=log,
xtick style={color=black},
y grid style={darkgray176},
ylabel={$ H_0^1\times L^2 $ error},
ymin=0.6e-11, ymax=10,
ymode=log,
ytick style={color=black}
]
\addplot [indianred204102119, mark=x, mark size=3, mark options={solid,fill opacity=0}, line width=1.2]
table {%
0.05 3.20757253135436
0.037037037037037 2.15834709134007
0.027027027027027 1.2718147767207
0.02 0.736962842933046
0.0147058823529411 0.416423206283391
0.010752688172043 0.232828161253726
0.0079365079365079 0.133813761521196
0.005813953488372 0.0781122816919549
0.0042918454935622 0.0493640549570749
0.0031545741324921 0.0350166520244034
0.0023201856148491 0.0287847010280547
0.0017064846416382 0.0263426522053124
0.0012562814070351 0.0254295296150551
0.0009242144177449 0.0250809684908338
0.000679809653297 0.0249412525331369
0.0005 0.0248812591298058
};
\addlegendentry{$ \text{CN} $}
\addplot [KITblue!80!black, mark=diamond, mark size=3, mark options={solid,fill opacity=0}, line width=1.5]
table {%
0.05 34617210539.6141
0.037037037037037 38144522682.5131
0.027027027027027 33415.5039674161
0.02 2383.88858006959
0.0147058823529411 17.5446094422665
0.010752688172043 0.230875963348323
0.0079365079365079 0.133534735272021
0.005813953488372 0.0780870057216154
0.0042918454935622 0.0493626603259111
0.0031545741324921 0.0350166277807897
0.0023201856148491 0.0287847137587837
0.0017064846416382 0.0263426664922899
0.0012562814070351 0.0254295444174681
0.0009242144177449 0.0250809834995358
0.000679809653297 0.0249412676275721
0.0005 0.0248812742600588
};
\addlegendentry{$ \text{DS}_{16}$}
\addplot [KITblue!80!black!35!white,dashed, mark=diamond*, mark size=3, mark options={solid,fill opacity=1.0}, line width=1.5, forget plot]
table {%
0.05 34617210539.612
0.037037037037037 38144522682.5255
0.027027027027027 33415.5039021574
0.02 2383.88841086831
0.0147058823529411 17.5408352603117
0.010752688172043 0.00623809471458192
0.0079365079365079 0.000361651593603734
0.005813953488372 3.26188301151212e-05
0.0042918454935622 2.01356330739032e-06
0.0031545741324921 6.30176424799033e-08
0.0023201856148491 8.88727542188057e-10
0.0017064846416382 8.71240645408515e-12
0.0012562814070351 4.8297177360466e-12
0.0009242144177449 1.56363650792585e-11
0.000679809653297 1.06047963495284e-11
0.0005 2.15077707129272e-11
};
\addplot [green!60!black, mark=*, mark size=3, mark options={solid,fill opacity=0}, line width=1.5]
table {%
0.05 1.20394478564752e+17
0.037037037037037 2.99606363614713e+17
0.027027027027027 5.07446291737767e+20
0.02 2.18197491403706e+22
0.0147058823529411 2.64444113902404e+17
0.010752688172043 5528855504.18053
0.0079365079365079 254448.861337818
0.005813953488372 0.31156505157938
0.0042918454935622 0.0493022227511072
0.0031545741324921 0.0350108953156318
0.0023201856148491 0.0287843943490013
0.0017064846416382 0.0263426563196961
0.0012562814070351 0.0254295442271116
0.0009242144177449 0.025080983497486
0.000679809653297 0.0249412676275694
0.0005 0.0248812742600811
};
\addlegendentry{$ \text{DS}_{8} $}
\addplot [green!60!black!35!white,dashed, mark=*, mark size=3, mark options={solid,fill opacity=1.0}, line width=1.5, forget plot]
table {%
				0.05 1.20394939724364e+17
				0.037037037037037 2.99607148858359e+17
				0.027027027027027 5.07447365892976e+20
				0.02 2.1819778774546e+22
				0.0147058823529411 2.64444283010434e+17
				0.010752688172043 5528857782.09616
				0.0079365079365079 254448.974172893
				0.005813953488372 0.301761846106904
				0.0042918454935622 0.000119249565436166
				0.0031545741324921 1.07285661782184e-05
				0.0023201856148491 8.51197810460674e-07
				0.0017064846416382 4.21476538124838e-08
				0.0012562814070351 1.24465333579348e-09
				0.0009242144177449 2.74811362740736e-11
				0.000679809653297 1.11883659444926e-11
				0.0005 2.1431743525803e-11
			};
\addplot [KITorange, mark=square, mark size=3, mark options={solid,fill opacity=0}, line width=1.5]
table {%
0.05 4.96355433269807e+29
0.037037037037037 2.51953378926535e+34
0.027027027027027 6.70162864067735e+36
0.02 1.87628843623342e+36
0.0147058823529411 5.0278985350251e+41
0.010752688172043 1.59887229097229e+47
0.0079365079365079 1.81574783765438e+44
0.005813953488372 3.52084058236585e+22
0.0042918454935622 1.55426669611737e+16
0.0031545741324921 31.272145333172
0.0023201856148491 0.0287753723837462
0.0017064846416382 0.0263417753706259
0.0012562814070351 0.0254294792320438
0.0009242144177449 0.0250809797079028
0.000679809653297 0.0249412674356959
0.0005 0.0248812742508294
}; 
\addlegendentry{$ \text{DS}_{4} $}
\addplot [KITorange50,dashed , mark=square*, mark size=3, mark options={solid,fill opacity=1.0}, line width=1.5, forget plot]
table {%
0.05 4.96355433269807e+29
0.037037037037037 2.51953378926535e+34
0.027027027027027 6.70162864067735e+36
0.02 1.87628843623342e+36
0.0147058823529411 5.0278985350251e+41
0.010752688172043 1.59887229097229e+47
0.0079365079365079 1.81574783765438e+44
0.005813953488372 3.52084058236585e+22
0.0042918454935622 1.55426669611737e+16
0.0031545741324921 31.2723984642368
0.0023201856148491 3.40041852983874e-05
0.0017064846416382 7.72582780755698e-06
0.0012562814070351 1.46648467579925e-06
0.0009242144177449 1.86678559282991e-07
0.000679809653297 1.52862524271891e-08
0.0005 8.49046814211984e-10
};
\addplot [black!80!white, mark=triangle, mark size=3, mark options={solid,rotate=180,fill opacity=0}, line width=1.2]
table {%
0.05 6.93165382358394e+74
0.037037037037037 5.04264132466292e+95
0.027027027027027 3.29177225018203e+125
0.02 inf
0.0147058823529411 inf
0.010752688172043 inf
0.0079365079365079 nan
0.005813953488372 nan
0.0042918454935622 nan
0.0031545741324921 nan
0.0023201856148491 nan
0.0017064846416382 nan
0.0012562814070351 nan
0.0009242144177449 1e99
0.000679809653297 0.0247882946094824
0.0005 0.0248029906815277
};
\addlegendentry{$ \text{lf} $}
\addplot [gray, dotted, line width=1.2]
table {%
0.05 1.0
0.0005 0.0001
};
\addlegendentry{$\landauO{\tau^{2}}$}
\end{axis}

\end{tikzpicture}
	\caption{Errors of the \ds\ approximation (DS$_\ell$)  with $ 4\times4 $ subdomains and $ \ell = 8 $ for the example from \Cref{subsec:2dex} measured at the final time $ T=1.0 $ against the exact solution (solid) and against the \CN\ approximation (dashed). For comparison, we also show the error of the
		\CN\ (CN) and the \LF\ (lf) methods on the full domain.}
	\label{fig::ex223single}
\end{figure}

\begin{figure}
	\begin{tikzpicture}

	\definecolor{brown1363485}{RGB}{136,34,85}
	\definecolor{burlywood221204119}{RGB}{221,204,119}
	\definecolor{cadetblue68170153}{RGB}{68,170,153}
	\definecolor{darkgray176}{RGB}{176,176,176}
	\definecolor{darkorchid17068153}{RGB}{170,68,153}
	\definecolor{darkslateblue5134136}{RGB}{51,34,136}
	\definecolor{forestgreen1711951}{RGB}{17,119,51}
	\definecolor{gainsboro221}{RGB}{221,221,221}
	\definecolor{indianred204102119}{RGB}{204,102,119}
	\definecolor{olivedrab15315351}{RGB}{153,153,51}
	\definecolor{skyblue136204238}{RGB}{136,204,238}

	\begin{axis}[
			width=0.96\textwidth,
			height=0.55\textwidth,
			legend cell align={left},
			legend style={
					fill opacity=0.8,
					draw opacity=1,
					text opacity=1,
					at={(0.97,0.03)},
					anchor=south east,
					draw=none
				},
			log basis x={10},
			log basis y={10},
			tick pos=both,
			unbounded coords=jump,
			x grid style={darkgray176},
			xlabel={\(\displaystyle \tau\)},
			xmin=0.000397164117362141, xmax=0.0629462705897084,
			xmode=log,
			xtick style={color=black},
			y grid style={darkgray176},
			ylabel={$H_0^1 \times L^2$ error},
			ymin=1e-11, ymax=10,
			ymode=log,
			ytick={1e0,1e-2,1e-5,1e-8,1e-10},
			ytick style={color=black}
		]
		\addplot [brown1363485, mark=o, mark size=3, line width=1.2, dashed,mark options={solid,fill opacity=0}]
		table {%
				0.05 223252659665158
				0.037037037037037 3.66968134362702e+17
				0.027027027027027 4.61864571500841e+20
				0.02 7.04406962175429e+21
				0.0147058823529411 209464553824838
				0.010752688172043 5751365010.5277
				0.0079365079365079 140322.400549056
				0.005813953488372 0.0169035835223383
				0.0042918454935622 5.35202306764424e-05
				0.0031545741324921 6.2655327047487e-06
				0.0023201856148491 5.03317670066139e-07
				0.0017064846416382 2.53438881828402e-08
				0.0012562814070351 7.71067116115139e-10
				0.0009242144177449 2.20979096887163e-11
				0.000679809653297 1.21401340666405e-11
				0.0005 2.48126633480102e-11
			};
		\addlegendentry{$\text{DS}_{8,2 \times 1}$}
		\addplot [burlywood221204119, mark=triangle, mark size=3,line width=1.2, dashed, mark options={solid,rotate=0,fill opacity=0}]
		table {%
				0.05 486108510093525
				0.037037037037037 8.79017317984869e+17
				0.027027027027027 1.35731364535294e+21
				0.02 3.00628408700159e+22
				0.0147058823529411 1.70575090098179e+15
				0.010752688172043 57906891144.6693
				0.0079365079365079 2556087.90216102
				0.005813953488372 0.795225950122827
				0.0042918454935622 8.97353438092479e-05
				0.0031545741324921 8.23652589578942e-06
				0.0023201856148491 6.53953673547418e-07
				0.0017064846416382 3.2432755420416e-08
				0.0012562814070351 9.66185842197438e-10
				0.0009242144177449 2.00953082266037e-11
				0.000679809653297 9.86387468681183e-12
				0.0005 1.33492045013777e-11
			};
		\addlegendentry{$\text{DS}_{8,8 \times 1}$}
		\addplot [skyblue136204238, mark=diamond*, mark size=4,line width=1.2, dashed, mark options={solid,rotate=180,fill opacity=0}]
		table {%
				0.05 2.27880349626716e+15
				0.037037037037037 6.36186423307021e+18
				0.027027027027027 1.55302545739189e+22
				0.02 4.82205857548146e+23
				0.0147058823529411 3.23975511119268e+16
				0.010752688172043 1602956819261.17
				0.0079365079365079 85695352.2819461
				0.005813953488372 21.5509734946058
				0.0042918454935622 0.00151577292162377
				0.0031545741324921 1.55895995079316e-05
				0.0023201856148491 1.22451565979281e-06
				0.0017064846416382 5.99343290350507e-08
				0.0012562814070351 1.7553439352335e-09
				0.0009242144177449 3.18944711417866e-11
				0.000679809653297 9.56828346672789e-12
				0.0005 1.25588457063978e-11
			};
		\addlegendentry{$\text{DS}_{8, 16 \times 1}$}
		\addplot [darkslateblue5134136, mark=Mercedes star, mark size=3.5,line width=1.2, dashed, mark options={solid,fill opacity=0}]
		table {%
				0.05 7.92388036087155e+15
				0.037037037037037 1.34835070918296e+19
				0.027027027027027 1.69478920552375e+22
				0.02 3.40391665128627e+23
				0.0147058823529411 1.13511463972888e+16
				0.010752688172043 298991543312.549
				0.0079365079365079 41542079.1116606
				0.005813953488372 10.2384897231229
				0.0042918454935622 0.0024149558776448
				0.0031545741324921 1.94204090607618e-05
				0.0023201856148491 1.52421914669643e-06
				0.0017064846416382 7.45047633006568e-08
				0.0012562814070351 2.17051358227538e-09
				0.0009242144177449 3.83670540575686e-11
				0.000679809653297 1.0105227084004e-11
				0.0005 1.23479949156288e-11
			};
		\addlegendentry{$\text{DS}_{8, 20 \times 1}$}
		\addplot [olivedrab15315351, mark=triangle, mark size=3,line width=1.2, dashed, mark options={solid,rotate=270,fill opacity=0}]
		table {%
				0.05 1.06932468538169e+17
				0.037037037037037 2.15609280739557e+17
				0.027027027027027 4.88585895627674e+20
				0.02 2.17899218444545e+22
				0.0147058823529411 2.64216194298768e+17
				0.010752688172043 5484854984.86231
				0.0079365079365079 252346.050450097
				0.005813953488372 0.301612186134582
				0.0042918454935622 8.3954857463499e-05
				0.0031545741324921 9.3888525971068e-06
				0.0023201856148491 7.50093629108422e-07
				0.0017064846416382 3.7539220493449e-08
				0.0012562814070351 1.12989603494271e-09
				0.0009242144177449 3.10209055931325e-11
				0.000679809653297 1.24189814889455e-11
				0.0005 3.00100843700095e-11
			};
		\addlegendentry{$\text{DS}_{8, 2 \times 2}$}
		\addplot [darkorchid17068153, mark=square, mark size=3,line width=1.2, dashed, mark options={solid,fill opacity=0}]
		table {%
				0.05 1.05989413225309e+17
				0.037037037037037 4.05059201681479e+18
				0.027027027027027 4.83437684136788e+21
				0.02 9.3157539121844e+22
				0.0147058823529411 6.96608732166932e+17
				0.010752688172043 135049331767.217
				0.0079365079365079 5435742.87854757
				0.005813953488372 1.66968329255035
				0.0042918454935622 0.00033839971843693
				0.0031545741324921 3.29606283744237e-06
				0.0023201856148491 2.42166288941624e-07
				0.0017064846416382 1.13867957954852e-08
				0.0012562814070351 3.16480173769322e-10
				0.0009242144177449 2.54811058583551e-11
				0.000679809653297 1.7276445159632e-11
				0.0005 3.67001135375201e-11
			};
		\addlegendentry{$\text{DS}_{8, 3 \times 3}$}
		\addplot [green!60!black, mark=*, mark size=2,line width=1.2, dashed, mark options={solid,fill opacity=0}]
		table {%
				0.05 1.20394939724364e+17
				0.037037037037037 2.99607148858359e+17
				0.027027027027027 5.07447365892976e+20
				0.02 2.1819778774546e+22
				0.0147058823529411 2.64444283010434e+17
				0.010752688172043 5528857782.09616
				0.0079365079365079 254448.974172893
				0.005813953488372 0.301761846106904
				0.0042918454935622 0.000119249565436166
				0.0031545741324921 1.07285661782184e-05
				0.0023201856148491 8.51197810460674e-07
				0.0017064846416382 4.21476538124838e-08
				0.0012562814070351 1.24465333579348e-09
				0.0009242144177449 2.74811362740736e-11
				0.000679809653297 1.11883659444926e-11
				0.0005 2.1431743525803e-11
			};
		\addlegendentry{$\text{DS}_{8, 4 \times 4}$}
		\addplot [indianred204102119, mark=x, mark size=3, mark options={solid,fill opacity=0},line width=1.2]
		table {%
				0.05 3.20757253147596
				0.037037037037037 2.15834709150891
				0.027027027027027 1.27181477700514
				0.02 0.736962843428935
				0.0147058823529411 0.416423207170886
				0.010752688172043 0.232828162853514
				0.0079365079365079 0.133813764317501
				0.005813953488372 0.0781122864953385
				0.0042918454935622 0.0493640625692333
				0.0031545741324921 0.0350166627646944
				0.0023201856148491 0.0287847140998464
				0.0017064846416382 0.0263426664926413
				0.0012562814070351 0.0254295444175375
				0.0009242144177449 0.0250809835002074
				0.000679809653297 0.0249412676272458
				0.0005 0.024881274260673
			};
		\addlegendentry{$\text{CN}$}
		\addplot [black!80!white, mark=triangle, mark size=3, mark options={solid,rotate=180,fill opacity=0}, line width=1.2]
		table {%
				0.05 6.9316538235954e+74
				0.037037037037037 5.04264132466485e+95
				0.027027027027027 3.29177225018256e+125
				0.02 inf
				0.0147058823529411 inf
				0.010752688172043 inf
				0.0079365079365079 nan
				0.005813953488372 nan
				0.0042918454935622 nan
				0.0031545741324921 nan
				0.0023201856148491 nan
				0.0017064846416382 nan
				0.0012562814070351 nan
				0.0009242144177449 1e99
				0.000679809653297 0.0247883097980372
				0.0005 0.0248030058608415
			};
		\addlegendentry{ $\text{lf}$}
	\end{axis}
\end{tikzpicture}
	\caption{Errors of the \ds\ approximation against the \CN\ approximation on different, rectangular subdomain configurations. All runs were made with overlap parameter $\ovp = 8$.}
	\label{fig::ex223topologies}
\end{figure}

We discretize $\fulldomain$ by a regular, triangular mesh containing $2\cdot 10^6$ cells with mesh width $ h = 10^{-3} $. The error is measured at the final time $T=1.0$. In \Cref{fig::ex223single} the convergence against the exact solution is depicted for the \LF scheme, the \CN method and the \ds method with different choices of the overlap parameter $\ovp$. Additionally, the norm of the difference between the final \ds approximation and the \CN solution is measured. The results are in line with the one-dimensional example in \Cref{Fig:1D_ell_dependence}.

In \Cref{fig::ex223topologies} compare the error of the \ds\ approximation against the \CN\ approximation on different, rectangular subdomain configurations comprising $N_x \times N_y$ subdomains in $x$ and $y$-direction, respectively, for a fixed overlap parameter $\ovp=8$. The approximations are denoted by DS$_{\ell,N_x \times N_y}$. While the error constants depend slightly on the number of neighboring subdomains, the stability is not affected by the number or configuration of subdomains as also stated by our error analysis.

\subsection{Experiments using a parallel implementation}

	In this subsection we present the performance of a parallel implementation of the domain splitting method. As in the previous subsections, we use FEniCSx \cite{Fenics} as the underlying framework. More specifically, we use the PETSc interface through petsc4py \cite{DalPetAl11}. The global mesh is generated using Gmsh \cite{GeuRe09} and then partitioned into non-overlapping subdomains using PT-Scotch \cite{CheP08}. 
	The parallel implementation relies on the construction of local meshes associated with the overlapping subdomains $\SDov{i}$ and the leapfrog prediction domains, as illustrated in \Cref{Fig:MeshesImplementation}. In addition, efficient local neighbor-to-neighbor communication between adjacent subdomains must be established. The implementation is consistently based on these local meshes $\Th[\SDov{i}]$, while direct access to the global mesh $\Th[\fulldomain]$ is avoided whenever possible.

	\begin{figure}
        \centering
        \begin{tikzpicture}[scale=0.5]

\coordinate (P0) at (0.0,0.0);
\coordinate (P1) at (8.0,0.0);
\coordinate (P2) at (8.0,5.0);
\coordinate (P3) at (0.0,5.0);
\coordinate (P4) at (0.7272727272713526,0.0);
\coordinate (P5) at (1.454545454542336,0.0);
\coordinate (P6) at (2.181818181812941,0.0);
\coordinate (P7) at (2.909090909083548,0.0);
\coordinate (P8) at (3.636363636354153,0.0);
\coordinate (P9) at (4.363636363626778,0.0);
\coordinate (P10) at (5.090909090901422,0.0);
\coordinate (P11) at (5.818181818176066,0.0);
\coordinate (P12) at (6.545454545450711,0.0);
\coordinate (P13) at (7.272727272725355,0.0);
\coordinate (P14) at (8.0,0.7142857142849608);
\coordinate (P15) at (8.0,1.428571428568591);
\coordinate (P16) at (8.0,2.142857142852804);
\coordinate (P17) at (8.0,2.857142857139607);
\coordinate (P18) at (8.0,3.571428571426405);
\coordinate (P19) at (8.0,4.285714285713203);
\coordinate (P20) at (7.272727272724246,5.0);
\coordinate (P21) at (6.545454545453539,5.0);
\coordinate (P22) at (5.818181818185859,5.0);
\coordinate (P23) at (5.09090909091818,5.0);
\coordinate (P24) at (4.3636363636505,5.0);
\coordinate (P25) at (3.636363636378783,5.0);
\coordinate (P26) at (2.909090909103026,5.0);
\coordinate (P27) at (2.181818181827269,5.0);
\coordinate (P28) at (1.454545454551513,5.0);
\coordinate (P29) at (0.7272727272757562,5.0);
\coordinate (P30) at (0.0,4.285714285714285);
\coordinate (P31) at (0.0,3.57142857142857);
\coordinate (P32) at (0.0,2.857142857142853);
\coordinate (P33) at (0.0,2.142857142857135);
\coordinate (P34) at (0.0,1.428571428571423);
\coordinate (P35) at (0.0,0.7142857142857117);
\coordinate (P36) at (0.6406971084116695,2.496786780131428);
\coordinate (P37) at (7.360322976148707,2.495589387098663);
\coordinate (P38) at (4.726161661606117,4.369587365236362);
\coordinate (P39) at (3.274726537514178,0.6287207136905271);
\coordinate (P40) at (2.544133927692886,4.370283672790174);
\coordinate (P41) at (5.454545454538743,0.6298366572994336);
\coordinate (P42) at (6.174140979179779,4.370689754642936);
\coordinate (P43) at (1.818181818177639,0.6298366572959363);
\coordinate (P44) at (6.867470873281542,0.6086417541664666);
\coordinate (P45) at (1.090909090913635,4.370163342699604);
\coordinate (P46) at (7.3494028615226,3.855735598708358);
\coordinate (P47) at (0.6547688164076052,1.139225694939238);
\coordinate (P48) at (1.806623297078301,4.369135070301849);
\coordinate (P49) at (2.165552901591474,3.74401452210516);
\coordinate (P50) at (2.903259317152124,3.74106067674268);
\coordinate (P51) at (2.531041954038659,3.118460716884476);
\coordinate (P52) at (3.264539607384775,3.112664202679856);
\coordinate (P53) at (2.900225990673177,2.489008913070915);
\coordinate (P54) at (3.62755704109661,2.48470550851066);
\coordinate (P55) at (2.167993293669921,2.492724346204462);
\coordinate (P56) at (3.991815802230595,3.110589057540085);
\coordinate (P57) at (4.353459784589963,2.484667396614063);
\coordinate (P58) at (4.717943890312323,3.109181614686769);
\coordinate (P59) at (5.079683227766391,2.479851588050829);
\coordinate (P60) at (5.443626727814732,3.106498214282779);
\coordinate (P61) at (5.80486699124863,2.478568030807522);
\coordinate (P62) at (4.716804475017599,1.862641115030466);
\coordinate (P63) at (2.539899902798354,1.86943739441172);
\coordinate (P64) at (6.178358572172582,3.106254332226279);
\coordinate (P65) at (1.812861992943986,1.867350784092799);
\coordinate (P66) at (6.555890540127068,2.482492574583282);
\coordinate (P67) at (6.173082117646791,1.857608065862945);
\coordinate (P68) at (6.887810002541,1.851267606102768);
\coordinate (P69) at (1.421631961522508,2.496508572865697);
\coordinate (P70) at (1.078321678387604,1.838228385733559);
\coordinate (P71) at (1.04518871783308,3.11692851310739);
\coordinate (P72) at (6.919242762810557,3.147726354150095);
\coordinate (P73) at (6.550056927047473,3.767467815863466);
\coordinate (P74) at (5.813892855468084,3.742435862118718);
\coordinate (P75) at (1.782257816415151,3.12387245431642);
\coordinate (P76) at (4.357101524241143,3.738564152152282);
\coordinate (P77) at (3.270487493896649,4.368748608779173);
\coordinate (P78) at (6.172024206958466,0.6199998102283633);
\coordinate (P79) at (5.45098581289999,4.369780147215225);
\coordinate (P80) at (5.447590848832858,1.858767811924908);
\coordinate (P81) at (2.546397198701698,0.6264179262653663);
\coordinate (P82) at (4.724112902303668,0.6242213325309243);
\coordinate (P83) at (3.996646733988828,0.6267715268962997);
\coordinate (P84) at (3.631405401422447,1.250869768280416);
\coordinate (P85) at (5.084952078723731,3.739341225948689);
\coordinate (P86) at (3.997221699792107,4.369583040829473);
\coordinate (P87) at (3.630862861924443,3.738521225859338);
\coordinate (P88) at (1.44452847870528,1.2215096578411);
\coordinate (P89) at (2.909396085873609,1.250536363762852);
\coordinate (P90) at (2.174027437031106,1.24657683510206);
\coordinate (P91) at (6.522374425257523,1.236656264708224);
\coordinate (P92) at (3.263178703252664,1.874386774249048);
\coordinate (P93) at (5.815040261984466,1.238261741191983);
\coordinate (P94) at (5.083749536477915,1.243663313938261);
\coordinate (P95) at (3.989208881955642,1.870588460054185);
\coordinate (P96) at (1.416582876566147,3.746484001266361);
\coordinate (P97) at (0.6504769182568116,3.824834745188991);
\coordinate (P98) at (4.353392356431248,1.246971245035119);
\coordinate (P99) at (6.894596821138778,4.412210941267588);
\coordinate (P100) at (1.103916492366339,0.5851999219475496);
\coordinate (P101) at (7.320836105287467,1.165185466531415);
\coordinate (P102) at (7.51063692170502,1.807842800367948);
\coordinate (P103) at (0.4747575206413757,1.809133886446557);
\coordinate (P104) at (7.524013956604464,3.214285714283007);
\coordinate (P105) at (0.4661104944371756,3.21428571428571);
\coordinate (P106) at (7.465452751953227,0.5311612121781986);
\coordinate (P107) at (0.5345472480476442,4.468838787821019);
\coordinate (P108) at (7.475798340905039,4.479372478390775);
\coordinate (P109) at (0.524201659093822,0.5206275216090231);

\definecolor{colorSD2}{RGB}{113,0,0}
\definecolor{colorSD3}{RGB}{0,0,113}

\fill[draw=colorSD3, fill=colorSD3,opacity = 0.8] 
(P0) -- (P35) -- (P34) -- (P33) -- (P33) -- (P36) -- (P69) -- (P55) -- (P53) -- (P54) -- (P57) -- (P62) -- (P94) -- (P82) -- (P9) -- (P8) -- (P7) -- (P6) -- (P5) -- (P4) -- cycle;

\fill[draw=colorSD3, fill=colorSD3,opacity = 0.5] 
(P33) -- (P36) -- (P69) -- (P55) -- (P53) -- (P54) -- (P57) -- (P62) -- (P94) -- (P82) -- (P9) -- (P10) -- (P11) -- (P78) -- (P91) -- (P67) -- (P61) -- (P60) -- (P85) -- (P76) -- (P87) -- (P50) -- (P49) -- (P96) -- (P97) -- (P31) -- cycle;

\foreach \a/\b/\c in {
            P72/P46/P73,P14/P15/P101,P70/P47/P88,P14/P101/P106,P46/P72/P104,P37/P66/P68,P47/P70/P103,P66/P37/P72,P69/P36/P70,P36/P69/P71,P73/P46/P99,P30/P31/P97,P88/P47/P100,P49/P48/P96,P55/P51/P75,P64/P66/P72,P64/P72/P73,P75/P49/P96,P69/P55/P75,P51/P49/P75,P37/P68/P102,P30/P97/P107,P73/P42/P74,P71/P69/P75,P48/P45/P96,P64/P73/P74,P71/P96/P97,P74/P42/P79,P81/P89/P90,P74/P79/P85,P76/P38/P86,P79/P38/P85,P12/P44/P78,P70/P36/P103,P40/P48/P49,P38/P76/P85,P25/P26/P77,P42/P22/P79,P40/P27/P48,P22/P23/P79,P82/P94/P98,P38/P24/P86,P64/P61/P66,P65/P55/P69,P65/P69/P70,P59/P62/P80,P60/P61/P64,P28/P29/P45,P12/P13/P44,P43/P6/P81,P28/P45/P48,P23/P38/P79,P66/P67/P68,P60/P59/P61,P41/P11/P78,P25/P77/P86,P63/P55/P65,P5/P6/P43,P21/P22/P42,P10/P41/P82,P26/P40/P77,P11/P12/P78,P53/P51/P55,P76/P86/P87,P26/P27/P40,P10/P11/P41,P86/P77/P87,P40/P49/P50,P24/P25/P86,P66/P61/P67,P27/P28/P48,P39/P83/P84,P7/P39/P81,P6/P7/P81,P58/P59/P60,P9/P10/P82,P39/P8/P83,P50/P49/P51,P58/P57/P59,P9/P82/P83,P61/P59/P80,P59/P57/P62,P53/P55/P63,P8/P9/P83,P7/P8/P39,P23/P24/P38,P60/P64/P74,P58/P60/P85,P56/P57/P58,P52/P53/P54,P89/P63/P90,P56/P54/P57,P52/P54/P56,P56/P58/P76,P50/P51/P52,P52/P51/P53,P76/P58/P85,P40/P50/P77,P60/P74/P85,P94/P62/P98,P50/P52/P87,P77/P50/P87,P67/P61/P80,P56/P76/P87,P52/P56/P87,P18/P19/P46,P34/P35/P47,P16/P17/P37,P32/P33/P36,P42/P73/P99,P81/P39/P89,P65/P70/P88,P78/P44/P91,P93/P80/P94,P43/P81/P90,P83/P82/P98,P53/P63/P92,P65/P88/P90,P62/P57/P95,P80/P62/P94,P67/P80/P93,P88/P43/P90,P39/P84/P89,P63/P65/P90,P62/P95/P98,P92/P84/P95,P63/P89/P92,P41/P93/P94,P20/P21/P99,P4/P5/P100,P68/P67/P91,P71/P75/P96,P54/P53/P92,P41/P78/P93,P82/P41/P94,P89/P84/P92,P54/P92/P95,P57/P54/P95,P101/P15/P102,P84/P83/P98,P96/P45/P97,P95/P84/P98,P72/P37/P104,P91/P67/P93,P43/P88/P100,P78/P91/P93,P36/P71/P105,P91/P44/P101,P21/P42/P99,P5/P43/P100,P13/P1/P106,P29/P3/P107,P1/P14/P106,P3/P30/P107,P68/P91/P101,P16/P37/P102,P34/P47/P103,P36/P33/P103,P18/P46/P104,P37/P17/P104,P32/P36/P105,P2/P20/P108,P0/P4/P109,P19/P2/P108,P35/P0/P109,P101/P44/P106,P71/P97/P105,P97/P45/P107,P15/P16/P102,P33/P34/P103,P97/P31/P105,P17/P18/P104,P31/P32/P105,P4/P100/P109,P20/P99/P108,P68/P101/P102,P44/P13/P106,P45/P29/P107,P46/P19/P108,P47/P35/P109,P99/P46/P108,P100/P47/P109}{
   \draw[thin,opacity=0.5] (\a) -- (\b) -- (\c) -- cycle;
}

\node[circle,draw=colorSD3,text=black,fill=white,fill opacity=0.9,inner sep = 0.1mm,scale=1.] at (P90) {$ \Omega_i $};
\draw[color=colorSD2,line width = 2.0] (P11) -- (P78) -- (P91) -- (P67) -- (P61) -- (P60) -- (P85) -- (P76) -- (P87) -- (P50) -- (P49) -- (P96) -- (P97) -- (P31) ; 
\node[circle,draw=colorSD3,text=black,fill=white,fill opacity=0.9,inner sep = 0.05mm,yshift=0.1mm,scale=1.] at (P51) {$ \Omega_i^{\delta} $};
\node[circle,draw=colorSD2,text=black,fill=white,fill opacity=0.9,inner sep = 0.1mm,scale=1.] at (P61) {$ \Gamma_i^{\ell} $};

\end{tikzpicture}
        \begin{tikzpicture}[scale=0.5]

\coordinate (P0) at (0.0,0.0);
\coordinate (P1) at (8.0,0.0);
\coordinate (P2) at (8.0,5.0);
\coordinate (P3) at (0.0,5.0);
\coordinate (P4) at (0.7272727272713526,0.0);
\coordinate (P5) at (1.454545454542336,0.0);
\coordinate (P6) at (2.181818181812941,0.0);
\coordinate (P7) at (2.909090909083548,0.0);
\coordinate (P8) at (3.636363636354153,0.0);
\coordinate (P9) at (4.363636363626778,0.0);
\coordinate (P10) at (5.090909090901422,0.0);
\coordinate (P11) at (5.818181818176066,0.0);
\coordinate (P12) at (6.545454545450711,0.0);
\coordinate (P13) at (7.272727272725355,0.0);
\coordinate (P14) at (8.0,0.7142857142849608);
\coordinate (P15) at (8.0,1.428571428568591);
\coordinate (P16) at (8.0,2.142857142852804);
\coordinate (P17) at (8.0,2.857142857139607);
\coordinate (P18) at (8.0,3.571428571426405);
\coordinate (P19) at (8.0,4.285714285713203);
\coordinate (P20) at (7.272727272724246,5.0);
\coordinate (P21) at (6.545454545453539,5.0);
\coordinate (P22) at (5.818181818185859,5.0);
\coordinate (P23) at (5.09090909091818,5.0);
\coordinate (P24) at (4.3636363636505,5.0);
\coordinate (P25) at (3.636363636378783,5.0);
\coordinate (P26) at (2.909090909103026,5.0);
\coordinate (P27) at (2.181818181827269,5.0);
\coordinate (P28) at (1.454545454551513,5.0);
\coordinate (P29) at (0.7272727272757562,5.0);
\coordinate (P30) at (0.0,4.285714285714285);
\coordinate (P31) at (0.0,3.57142857142857);
\coordinate (P32) at (0.0,2.857142857142853);
\coordinate (P33) at (0.0,2.142857142857135);
\coordinate (P34) at (0.0,1.428571428571423);
\coordinate (P35) at (0.0,0.7142857142857117);
\coordinate (P36) at (0.6406971084116695,2.496786780131428);
\coordinate (P37) at (7.360322976148707,2.495589387098663);
\coordinate (P38) at (4.726161661606117,4.369587365236362);
\coordinate (P39) at (3.274726537514178,0.6287207136905271);
\coordinate (P40) at (2.544133927692886,4.370283672790174);
\coordinate (P41) at (5.454545454538743,0.6298366572994336);
\coordinate (P42) at (6.174140979179779,4.370689754642936);
\coordinate (P43) at (1.818181818177639,0.6298366572959363);
\coordinate (P44) at (6.867470873281542,0.6086417541664666);
\coordinate (P45) at (1.090909090913635,4.370163342699604);
\coordinate (P46) at (7.3494028615226,3.855735598708358);
\coordinate (P47) at (0.6547688164076052,1.139225694939238);
\coordinate (P48) at (1.806623297078301,4.369135070301849);
\coordinate (P49) at (2.165552901591474,3.74401452210516);
\coordinate (P50) at (2.903259317152124,3.74106067674268);
\coordinate (P51) at (2.531041954038659,3.118460716884476);
\coordinate (P52) at (3.264539607384775,3.112664202679856);
\coordinate (P53) at (2.900225990673177,2.489008913070915);
\coordinate (P54) at (3.62755704109661,2.48470550851066);
\coordinate (P55) at (2.167993293669921,2.492724346204462);
\coordinate (P56) at (3.991815802230595,3.110589057540085);
\coordinate (P57) at (4.353459784589963,2.484667396614063);
\coordinate (P58) at (4.717943890312323,3.109181614686769);
\coordinate (P59) at (5.079683227766391,2.479851588050829);
\coordinate (P60) at (5.443626727814732,3.106498214282779);
\coordinate (P61) at (5.80486699124863,2.478568030807522);
\coordinate (P62) at (4.716804475017599,1.862641115030466);
\coordinate (P63) at (2.539899902798354,1.86943739441172);
\coordinate (P64) at (6.178358572172582,3.106254332226279);
\coordinate (P65) at (1.812861992943986,1.867350784092799);
\coordinate (P66) at (6.555890540127068,2.482492574583282);
\coordinate (P67) at (6.173082117646791,1.857608065862945);
\coordinate (P68) at (6.887810002541,1.851267606102768);
\coordinate (P69) at (1.421631961522508,2.496508572865697);
\coordinate (P70) at (1.078321678387604,1.838228385733559);
\coordinate (P71) at (1.04518871783308,3.11692851310739);
\coordinate (P72) at (6.919242762810557,3.147726354150095);
\coordinate (P73) at (6.550056927047473,3.767467815863466);
\coordinate (P74) at (5.813892855468084,3.742435862118718);
\coordinate (P75) at (1.782257816415151,3.12387245431642);
\coordinate (P76) at (4.357101524241143,3.738564152152282);
\coordinate (P77) at (3.270487493896649,4.368748608779173);
\coordinate (P78) at (6.172024206958466,0.6199998102283633);
\coordinate (P79) at (5.45098581289999,4.369780147215225);
\coordinate (P80) at (5.447590848832858,1.858767811924908);
\coordinate (P81) at (2.546397198701698,0.6264179262653663);
\coordinate (P82) at (4.724112902303668,0.6242213325309243);
\coordinate (P83) at (3.996646733988828,0.6267715268962997);
\coordinate (P84) at (3.631405401422447,1.250869768280416);
\coordinate (P85) at (5.084952078723731,3.739341225948689);
\coordinate (P86) at (3.997221699792107,4.369583040829473);
\coordinate (P87) at (3.630862861924443,3.738521225859338);
\coordinate (P88) at (1.44452847870528,1.2215096578411);
\coordinate (P89) at (2.909396085873609,1.250536363762852);
\coordinate (P90) at (2.174027437031106,1.24657683510206);
\coordinate (P91) at (6.522374425257523,1.236656264708224);
\coordinate (P92) at (3.263178703252664,1.874386774249048);
\coordinate (P93) at (5.815040261984466,1.238261741191983);
\coordinate (P94) at (5.083749536477915,1.243663313938261);
\coordinate (P95) at (3.989208881955642,1.870588460054185);
\coordinate (P96) at (1.416582876566147,3.746484001266361);
\coordinate (P97) at (0.6504769182568116,3.824834745188991);
\coordinate (P98) at (4.353392356431248,1.246971245035119);
\coordinate (P99) at (6.894596821138778,4.412210941267588);
\coordinate (P100) at (1.103916492366339,0.5851999219475496);
\coordinate (P101) at (7.320836105287467,1.165185466531415);
\coordinate (P102) at (7.51063692170502,1.807842800367948);
\coordinate (P103) at (0.4747575206413757,1.809133886446557);
\coordinate (P104) at (7.524013956604464,3.214285714283007);
\coordinate (P105) at (0.4661104944371756,3.21428571428571);
\coordinate (P106) at (7.465452751953227,0.5311612121781986);
\coordinate (P107) at (0.5345472480476442,4.468838787821019);
\coordinate (P108) at (7.475798340905039,4.479372478390775);
\coordinate (P109) at (0.524201659093822,0.5206275216090231);

\definecolor{colorSD2}{RGB}{113,0,0}
\definecolor{colorSD3}{RGB}{0,113,0}

\fill[draw=colorSD3, fill=colorSD3,opacity = 0.5] 
(P31) -- (P32) -- (P105) -- (P71) -- (P75) -- (P51) -- (P52) -- (P56) -- (P58) -- (P59) -- (P80) -- (P93) -- (P41) -- (P10) -- (P11) -- (P12) -- (P44) -- (P101) -- (P68) -- (P66) -- (P64) -- (P74) -- (P79) -- (P38) -- (P86) -- (P77) -- (P40) -- (P48) -- (P45) -- (P107) -- (P30) -- cycle;

\foreach \a/\b/\c in {
            P72/P46/P73,P14/P15/P101,P70/P47/P88,P14/P101/P106,P46/P72/P104,P37/P66/P68,P47/P70/P103,P66/P37/P72,P69/P36/P70,P36/P69/P71,P73/P46/P99,P30/P31/P97,P88/P47/P100,P49/P48/P96,P55/P51/P75,P64/P66/P72,P64/P72/P73,P75/P49/P96,P69/P55/P75,P51/P49/P75,P37/P68/P102,P30/P97/P107,P73/P42/P74,P71/P69/P75,P48/P45/P96,P64/P73/P74,P71/P96/P97,P74/P42/P79,P81/P89/P90,P74/P79/P85,P76/P38/P86,P79/P38/P85,P12/P44/P78,P70/P36/P103,P40/P48/P49,P38/P76/P85,P25/P26/P77,P42/P22/P79,P40/P27/P48,P22/P23/P79,P82/P94/P98,P38/P24/P86,P64/P61/P66,P65/P55/P69,P65/P69/P70,P59/P62/P80,P60/P61/P64,P28/P29/P45,P12/P13/P44,P43/P6/P81,P28/P45/P48,P23/P38/P79,P66/P67/P68,P60/P59/P61,P41/P11/P78,P25/P77/P86,P63/P55/P65,P5/P6/P43,P21/P22/P42,P10/P41/P82,P26/P40/P77,P11/P12/P78,P53/P51/P55,P76/P86/P87,P26/P27/P40,P10/P11/P41,P86/P77/P87,P40/P49/P50,P24/P25/P86,P66/P61/P67,P27/P28/P48,P39/P83/P84,P7/P39/P81,P6/P7/P81,P58/P59/P60,P9/P10/P82,P39/P8/P83,P50/P49/P51,P58/P57/P59,P9/P82/P83,P61/P59/P80,P59/P57/P62,P53/P55/P63,P8/P9/P83,P7/P8/P39,P23/P24/P38,P60/P64/P74,P58/P60/P85,P56/P57/P58,P52/P53/P54,P89/P63/P90,P56/P54/P57,P52/P54/P56,P56/P58/P76,P50/P51/P52,P52/P51/P53,P76/P58/P85,P40/P50/P77,P60/P74/P85,P94/P62/P98,P50/P52/P87,P77/P50/P87,P67/P61/P80,P56/P76/P87,P52/P56/P87,P18/P19/P46,P34/P35/P47,P16/P17/P37,P32/P33/P36,P42/P73/P99,P81/P39/P89,P65/P70/P88,P78/P44/P91,P93/P80/P94,P43/P81/P90,P83/P82/P98,P53/P63/P92,P65/P88/P90,P62/P57/P95,P80/P62/P94,P67/P80/P93,P88/P43/P90,P39/P84/P89,P63/P65/P90,P62/P95/P98,P92/P84/P95,P63/P89/P92,P41/P93/P94,P20/P21/P99,P4/P5/P100,P68/P67/P91,P71/P75/P96,P54/P53/P92,P41/P78/P93,P82/P41/P94,P89/P84/P92,P54/P92/P95,P57/P54/P95,P101/P15/P102,P84/P83/P98,P96/P45/P97,P95/P84/P98,P72/P37/P104,P91/P67/P93,P43/P88/P100,P78/P91/P93,P36/P71/P105,P91/P44/P101,P21/P42/P99,P5/P43/P100,P13/P1/P106,P29/P3/P107,P1/P14/P106,P3/P30/P107,P68/P91/P101,P16/P37/P102,P34/P47/P103,P36/P33/P103,P18/P46/P104,P37/P17/P104,P32/P36/P105,P2/P20/P108,P0/P4/P109,P19/P2/P108,P35/P0/P109,P101/P44/P106,P71/P97/P105,P97/P45/P107,P15/P16/P102,P33/P34/P103,P97/P31/P105,P17/P18/P104,P31/P32/P105,P4/P100/P109,P20/P99/P108,P68/P101/P102,P44/P13/P106,P45/P29/P107,P46/P19/P108,P47/P35/P109,P99/P46/P108,P100/P47/P109}{
   \draw[thin,opacity=0.5] (\a) -- (\b) -- (\c) -- cycle;
}

\draw[color=colorSD2,line width = 2.0] (P11) -- (P78) -- (P91) -- (P67) -- (P61) -- (P60) -- (P85) -- (P76) -- (P87) -- (P50) -- (P49) -- (P96) -- (P97) -- (P31) ; 
\node[circle,draw=colorSD2,text=black,fill=white,fill opacity=0.9,inner sep = 0.1mm,scale=1.] at (P61) {$ \Gamma_i^{\ell} $};

\end{tikzpicture}
        \caption{Overlapping subdomain $\SDov{i}$ (left, dark and light blue area) for $\ovp=2$ and leapfrog prediction domain (right, green area).
        The interface $\SDovint{i}$ is shown in dark red in both figures.
        Only the colored parts of the mesh are stored locally.
        }
        \label{Fig:MeshesImplementation}
    \end{figure}

\begin{figure}[t]
		\includegraphics[width=0.4\textwidth]{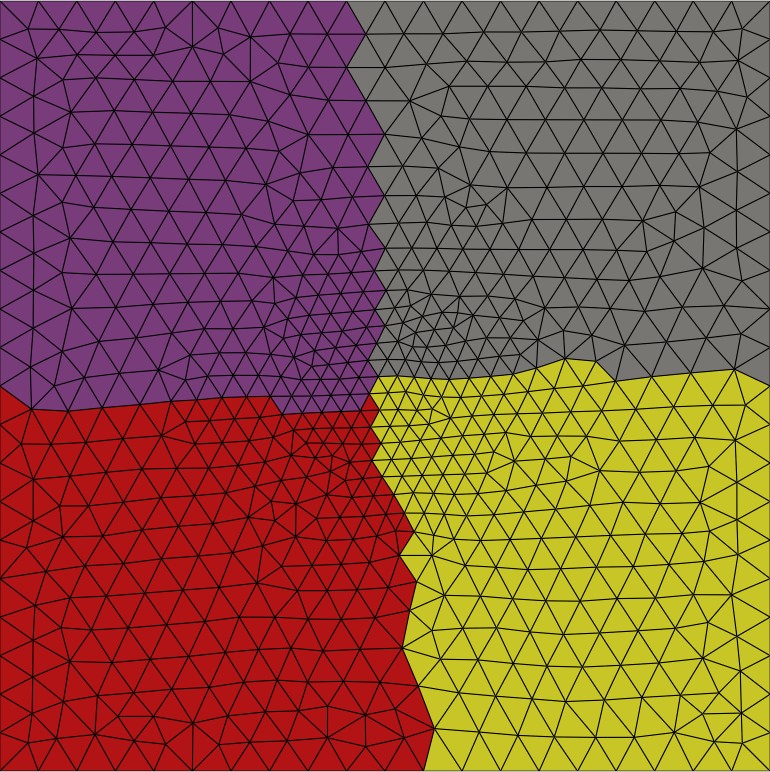}
		\caption{Partitioning of the global domain into four nonoverlapping subdomains using PT-scotch with $h_{\min}=0.0215$. For the simulations, significantly finer meshes were used.}
		\label{Fig:ScotchPartitioning}
	\end{figure}

	The global mesh is constructed on the unit square $[0,1]^2$, which we refine around the center with a refinement factor of $2$. An example of such a mesh, partitioned into four non-overlapping subdomains using PT-Scotch, is shown in \Cref{Fig:ScotchPartitioning}.

	We consider a global mesh with $h_{\min} \approx 2.9 \cdot 10^{-4}$, resulting in approximately $4.89 \cdot 10^{6}$ spatial degrees of freedom when accounting for both components.
	We compare the domain splitting method with $\ovp$ layers ($\text{DS}_{\ovp}$) with parallelized global time-integration schemes, namely the \CN (CN) and the leapfrog (LF) methods, which exploit the built-in parallelization provided by FEniCSx. For the  \ds and the \CN method, the resulting linear systems are solved using the PETSc Cholesky decomposition of the system matrix, which proved to be faster for the problem size considered here than the available iterative solvers provided by petsc4py. The implementation of the leapfrog method and the prediction step within DS are fully explicit, since the mass matrix is diagonal.

	To evaluate the accuracy of the different approaches, we prescribe an exact solution and adapt the right-hand side accordingly. Specifically, we choose
	\begin{equation*}
			u(x,y,t) = \sin(\pi x)^2 \sin(\pi y)^2 \exp(t), \quad t \in [0,5],
	\end{equation*}
	which allows us to compute the relative $L^2$ error of the $u$ component of a given approximation compared to the $L^2$ norm of the exact solution. All simulations are performed on a single node using eight CPU cores of an Intel Xeon Platinum 8358 processor. We show the relative $L^2$ errors in \Cref{Fig:ParallelConvergence} and compare the simulation timings in \Cref{Table:comparisonDSCN}.

	\begin{figure}
		\begin{tikzpicture}

\definecolor{burlywood221204119}{RGB}{221,204,119}
\definecolor{darkgray176}{RGB}{176,176,176}
\definecolor{darkslateblue5134136}{RGB}{51,34,136}
\definecolor{forestgreen1711951}{RGB}{17,119,51}
\definecolor{gray}{RGB}{128,128,128}
\definecolor{indianred204102119}{RGB}{204,102,119}

\begin{axis}[
	width=0.96\textwidth,
	height=0.4\textwidth,
legend cell align={left},
legend columns=3,
legend style={
  fill opacity=0.0,
  draw opacity=1,
  text opacity=1,
  at={(0.97,0.97)},
  anchor=north east,
  draw=none
},
log basis x={10},
log basis y={10},
tick pos=both,
unbounded coords=jump,
x grid style={darkgray176},
xlabel={\(\displaystyle \tau\)},
xmin=0.00008, xmax=0.12,
xmode=log,
xtick style={color=black},
y grid style={darkgray176},
ylabel={relative $ L^2 $ error in $u$},
ymin=5e-7, ymax=1e-3,
ymode=log,
ytick style={color=black}
]
\addplot [indianred204102119, mark=x, mark size=3, mark options={solid,fill opacity=0}, line width=1.2]
table {%
  0.1 8.005099976892411e-05
  0.046296296296296294 1.882005783407872e-05
  0.021739130434782608 5.841532395592435e-06
  0.01 3.096390636210233e-06
  0.004651162790697674 2.5299075071123466e-06
  0.002154243860404998 2.410486956970758e-06
  0.001 2.3849397308136913e-06
  0.0004641663572224285 2.3794354263869967e-06
  0.00021544295070665288 2.378217346206064e-06
  0.0001 2.3777862674072384e-06
};
\addlegendentry{$ \text{CN} $}

\addplot [green!60!black, mark=*, mark size=3.5, mark options={solid,fill opacity=0}, line width=1.0]
table {%
    0.004642525533890436 1.3116933818380909e+137
    0.002154243860404998 4.1802659505360544e-05
    0.001 2.3528362174576554e-06
    0.0004641663572224285 2.3794164266935185e-06
    0.00021544295070665288 2.3781774539375386e-06
    0.0001 2.3778049977241862e-06
};
\addlegendentry{$ \text{DS}_{16} $}

\addplot [blue!80!black, mark=diamond, mark size=3.5, mark options={solid,fill opacity=0}, line width=1.0]
table {%
    0.002154243860404998 1.5826691987895858e+77
    0.001 1.4598097543921503e-05
    0.0004641663572224285 2.359931394141873e-06
    0.00021544295070665288 2.378176468112019e-06
    0.0001 2.3778022349318347e-06
};
\addlegendentry{$ \text{DS}_{8}$}
\addplot [orange, mark=square, mark size=2, mark options={solid,fill opacity=0}, line width=1.0]
table {%
    0.001 3.5934888131137365e+57
    0.0004641663572224285 8.840755269914931e-06
    0.00021544295070665288 2.318366111912583e-06
    0.0001 2.377615663126389e-06
}; 
\addlegendentry{$ \text{DS}_{4} $}
\addplot [black!80!white, mark=triangle, mark size=3, mark options={solid,rotate=180,fill opacity=0}, line width=1.2]
		table {%
        0.00021542438604049978 1e100
        0.0001 2.379937162001152e-06
        };
\addlegendentry{ $\text{lf}$}
\addplot [gray, dotted, line width=1.2]
table {%
0.1 0.5e-4
0.0001 0.5e-10
};
\addlegendentry{$\mathcal{O}(\tau^{2})$}
\end{axis}

\end{tikzpicture}
 
		\caption{Comparison of the relative $L^2$ error in the first component for the domain splitting approximation ($\text{DS}_{\ovp}$), and the approximation obtained by the global \CN method (CN) and the leapfrog method (lf).}
		\label{Fig:ParallelConvergence}
	\end{figure}

	\setlength{\tabcolsep}{6pt}
     \begin{table}[H]
        \centering
        \begin{tabular}[t]{l|rrrrr}
        \hline
        &\multicolumn{1}{c}{$\text{lf}$} 
		& \multicolumn{1}{c}{$\text{DS}_4$} 
		& \multicolumn{1}{c}{$\text{DS}_8$} 
		& \multicolumn{1}{c}{$\text{DS}_{16}$} 
		& \multicolumn{1}{c}{CN}\\
        \hline
        $\tau = 1.00 \cdot 10^{-4}$ & 7\,367.0 s & 11\,525.5 s & 11\,186.2 s & 12\,026.6 s & 17\,267.1 s \\
        $\tau \approx 2.15\cdot 10^{-4}$ &\multicolumn{1}{c}{--} &5\,169.2 s &  5\,291.7 s & 5\,488.7 s & 7\,338.1 s \\
        $\tau \approx 4.64 \cdot 10^{-4}$ &\multicolumn{1}{c}{--} &\multicolumn{1}{c}{--} & 2\,473.0 s & 2\,294.8 s & 3\,363.1 s \\
        $\tau = 1.00 \cdot 10^{-3 }$ &\multicolumn{1}{c}{--} &\multicolumn{1}{c}{--} & \multicolumn{1}{c}{--} & 1\,259.9 s & 1\,380.3 s \\
        \hline
        \end{tabular}
        \caption{Comparison of timings (setup time for solver $+$ wall time for time loop) for the simulations from \Cref{Fig:ParallelConvergence}. We only included simulations, where the error is comparable to the error obtained from the global \CN method. Due to the dominance of the space discretization error, all corresponding relative $L^2$ errors lie between $2.35\cdot 10^{-6}$ and $2.38\cdot 10^{-6}$.}
        \label{Table:comparisonDSCN}
    \end{table}%

	For DS with eight subdomains, $\ell = 16$, and a time step $\tau = 0.001$, we obtain a relative $L^2$ error of $2.35\cdot 10^{-6}$ with a total runtime of $1\,259.9$ seconds. 
	In comparison, the parallel CN method using eight cores and the same time step yields a relative $L^2$ error of $2.38\cdot 10^{-6}$ with a runtime of $1\,380.3$ seconds.
	In general, for a fixed \timestepsize $\tau$, if the domain splitting scheme omits stability, it appears to be faster than the parallelized global CN method. 
	Moreover, an increase in the number of layers $\ovp$ within the overlap comes with an additional overhead, mainly due to communication, which is visible in most cases by an increased total simulation time. 
	Hence, we recommend to use $\ovp$ as small as possible while still ensuring stability.
	Finally, the parallel leapfrog scheme with eight cores and a smaller time step $\tau = 0.0001$ achieves a comparable relative $L^2$ error of $2.37\cdot10^{-6}$ at a total runtime of $7\,367.0$ seconds.

	The experiments indicate that the parallel implementation of the \ds method achieves accuracy
	comparable to the global methods, while improving stability compared to explicit methods and reducing computational costs compared to implicit methods.

\appendix
\section{Bounds of Theorem \texorpdfstring{\ref{ExponentialDecay}}{}}
\label{sec::Appendix}
	
	With the notation introduced in \Cref{def::bilinearforms} and in the proof of \Cref{ExponentialDecay} we now prove two Lemmas  used there. Recall the definition of the constants $\CML$ and $\Cint$ from \eqref{const::CML} and \eqref{const::CIh}, respectively, and the definition of the weight function $\weight$ in \eqref{eq::DefWeightFunction}.
	Let $\abs{\cdot}_{2,\K, \weight(-\Ipd) }$ denote the $H^2$ semi-norm, weighted by $\weight(-\Ipd)$, in the same manner as the bilinear forms from \Cref{def::bilinearforms}.
	Moreover, $\abs{\cdot}$ denotes the Euclidean norm when applied to a  vector.
	
	\begin{lemma}
		\label{Lem::BoundDecayTerm2}
		Let $\gamma \in (0,1]$, $\lambda > 0$, $ z \in \Vh[\SDov{i}]$. For $\wz = \weight(\Ipd)z$ defined in \eqref{def::wz} the interpolation error $ \IpdwzE = \wz - \NodalInt \wz$ satisfies
		\begin{equation*}
			\abs{\ba{\SDov{i}}\ip{\IpdwzE,z}} 
			\leq \gamma C_1 \bbw{\SDov{i}}{\weight(\Ipd) } \ip{z,z},
		\end{equation*}
		with $ C_1 = \ps \Cint \sqrt{ 3 \dim \e \cdot \max \{ \e\CML , 2\}}$.
	\end{lemma}
	\begin{proof}
	
	By definition of $\baSDov{i}\ip{\cdot,\cdot}$ in \eqref{def::ba}, $\baw{\SDov{i}}{\weight}(\cdot,\cdot)$ in \eqref{def::baw}, the \CSI, and \eqref{eq::weightProps-b} we obtain
	\begin{equation}
		\label{eq::weightedCSI}
		\begin{aligned}
			\abs{\ba{\SDov{i}}\ip{\IpdwzE,z}} 
			&= \abs{\ipD{\SDov{i}}{\ps^2 \grad (\IpdwzE), \grad z}} \\
			&= \abs{\ip{\ps^2 \weight(-\Ipd) ^{\half} \grad (\IpdwzE), \weight(\Ipd) ^{\half} \grad z}_{\SDov{i}}} \\
			&\leq \baw{\SDov{i}}{\weight(-\Ipd) }\ip{\IpdwzE, \IpdwzE}^{\half} \baw{\SDov{i}}{\weight(\Ipd) }\ip{z,z}^{\half}.
		\end{aligned}
	\end{equation}
	We estimate the first factor via
	\begin{align}
		\nonumber \baw{\SDov{i}}{\weight(-\Ipd) }\ip{\IpdwzE, \IpdwzE} 
			&= \sumTh[\SDov{i}] \baw{\K}{\weight(-\Ipd) }\ip{\IpdwzE, \IpdwzE} \\
			&\leq \sumTh[\SDov{i}] \max_{x \in \K} \weight(-\Ipd)  \ps^2 \abs{\IpdwzE}_{1,\K}^2 .
			\label{eq::goDownOnCells}
	\end{align}
	Since $\restr{z}{\K} \in P_1$, for $\K \in \Th[\SDov{i}]$ the function $ \restr{\wz}{\K} $ is smooth, so that the interpolation estimate $\eqref{const::CIh}$  yields
	\begin{equation}
		\label{eq::interpolationEst-weight}
		\abs{\IpdwzE}_{1,\K}^2 
		\leq \Cint^2 h^2 \abs{\wz}_{2,\K}^2 \leq \frac{\Cint^2 h^2}{\min_{x \in \K} \weight(-\Ipd)} \abs{\wz}_{2,\K, \weight(-\Ipd) }^2.
	\end{equation}
	Moreover, since $\weight(x)$ is monotonically increasing, there holds on $K\in \Th[\SDov{i}]$
	\begin{equation}
		\label{eq::weightMaxByMin}
		\begin{aligned}
			\frac{\max_{x \in \K} \weight(-\Ipd) }{\min_{x \in \K} \weight(-\Ipd) } 
			= \frac{\weight(- \min_{x \in \K} \Ipd)}{\weight( - \max_{x \in \K} \Ipd)}
			= \weight(\max_{x \in \K} \Ipd - \min_{x \in \K} \Ipd)  
			 \leq \weight(\diam \K) 
			\leq \weight(h) .
		\end{aligned}
	\end{equation}
	By combining this with \eqref{eq::goDownOnCells} and \eqref{eq::interpolationEst-weight} we get
	\begin{equation}
		\label{b.4}
		\begin{aligned}
			\baw{\SDov{i}}{\weight(-\Ipd) }\ip{\IpdwzE, \IpdwzE} 
			&\leq \ps^2 \Cint^2 h^2 \weight(h) \sumTh[\SDov{i}] \abs{\wz}_{2,\K, \weight(-\Ipd) }^2 .
		\end{aligned}
	\end{equation}
	Since $\d{r}(\d{j} \restr{v_h}{K}) = 0$ for $v_h \in P_1$ and $r,j = 1, \dots , \dim$, we get 
	\begin{equation}
		\label{eq::differentiated2wz}
		\frac{\partial^2}{\partial_r \partial_j} \wz = \frac{\gamma}{\max{\lambda, h}} \Bigl( \d{r} z \d{j} \Ipd + z \frac{\gamma}{\max\{ \lambda, h \} } \d{r}\Ipd \d{j} \Ipd  + \d{j} z \d{r} \Ipd \Bigr)\weight(\Ipd).
	\end{equation}
	Note, that $d: x \mapsto \dist(x, \SDovint{i})$ satisfies
	\begin{equation*}
		\abs{\grad d(\widehat{x})} = 1, \quad  \text{for all } \widehat{x} \in \SDov{i},
	\end{equation*}
	cf. \cite{CafCr10}. By definition of $d_h$ in \eqref{eq::DefInterpolatedDistance}, we deduce with a geometric argument
	\begin{equation}
		\label{eq::Derivdh}
		\abs{\partial_j d_h(\widehat{x})} \leq \abs{\grad d_h(\widehat{x})} \leq \abs{\grad d(\widehat{x})} = 1, \quad \text{ for } j=1,\dots,m
	\end{equation}
	 for all $\widehat{x} \in \SDov{i}$, see also \cite{BluLR92} after equation (24).
	
	By combining \eqref{eq::differentiated2wz} and \eqref{eq::Derivdh} and using the linearity of the integral we have
	\begin{equation}
		\label{b.5}
		\begin{aligned}
		\sumTh[\SDov{i}] &\abs{\wz}_{2,\K,\weight(-\Ipd) }^2 
		= \sumTh[\SDov{i}] \int_{\K} \weight(-\Ipd)  \sum_{r,j=1}^m \abs{\frac{\partial^2 \wz}{\partial_r \partial_j}}^2 \dintx \\
		&\leq \frac{ 3 \dim \gamma^2}{\max \{ \lambda, h\}^4} \left( \gamma^2 \ipD{\SDov{i}}{\weight(\Ipd)  z, z} + 2 \max \{ \lambda, h\}^2 \baw{\SDov{i}}{\weight(\Ipd) } \ip{z,z} \right) .
		\end{aligned}
	\end{equation}
	Using \eqref{const::CML}, one can easily show, that
	\begin{equation}
		\label{b.6}
		\begin{aligned}
			\ipSDov{i}{\weight(\Ipd) u,u} 
			\leq \CML \weight(h) \ipML{\SDov{i}}{\weight(\Ipd) u,u}, \qquad u \in \Vh[\SDov{i}].
		\end{aligned}
	\end{equation}
	Combining this with \eqref{b.5} yields 
	\begin{equation}
		\label{b.7}
		\begin{aligned}
			\sumTh[\SDov{i}] \abs{\wz}_{2,\K,\weight(-\Ipd) }^2 
			&\leq \frac{3 \dim \gamma^2}{\max \{ \lambda, h\}^2} \max \{ \gamma^2 \CML \weight(h), 2\} \bbw{\SDov{i}}{\weight(\Ipd) } \ip{z,z} .
		\end{aligned}
	\end{equation}
	Putting \eqref{b.4}, \eqref{b.7} together yields an estimate for the first factor of \eqref{eq::weightedCSI} of the form
	\begin{equation*}
		\baw{\SDov{i}}{\weight(-\Ipd) }\ip{\IpdwzE, \IpdwzE} \leq \ps^2 \Cint^2 h^2 \weight(h) \frac{3 \dim \gamma^2}{\max \{ \lambda, h\}^2} \max \{ \gamma^2 \CML \weight(h), 2\} \bbw{\SDov{i}}{\weight(\Ipd) } \ip{z,z}.
	\end{equation*}
	By \eqref{def::bbw} and since $ \frac{1}{\lambda^2}\ipML{\SDov{i}}{\weight(\Ipd) z,z} \geq 0 $, we bound the second factor in \eqref{eq::weightedCSI} by
	\begin{equation*}
		\baw{\SDov{i}}{\weight(\Ipd) }\ip{z,z} \leq \bbw{\SDov{i}}{\weight(\Ipd) } \ip{z,z}.
	\end{equation*}
	In total, we receive the estimate
	\begin{equation*}
		\abs{\ba{\SDov{i}}\ip{\IpdwzE,z}} 
		\leq \ps \Cint \weight(\frac{h}{2}) \frac{\sqrt{3 \dim }\gamma h }{\max \{\lambda, h\}} \max \{ \gamma \CML^{\half} \weight(\frac{h}2), \sqrt{2}\} \bbw{\SDov{i}}{\weight(\Ipd) } \ip{z,z}.
	\end{equation*}
	The claim follows by using that $ \frac{h}{\max \{ \lambda, h \}} \leq 1 $, $\gamma \leq 1$, and thus $ \weight(h/2 ) \leq \e^{\gamma/2} \leq \sqrt{\e} $.
	\end{proof}

	By similar arguments as in the proof of \Cref{Lem::BoundDecayTerm2} we can also derive a bound for the last remaining term in \eqref{eq::weightedCombinedNorm}.
	\begin{lemma}
		\label{Lem::BoundDecayTerm3}
		Let the assumptions of \Cref{Lem::BoundDecayTerm2} hold.
		Then, we have 
		\begin{equation*}
			\abs{\int_{\SDov{i}} z \ps^2 \grad \weight(\Ipd) \cdot \grad z \dintx  }  \leq \gamma C_2 \bbw{\SDov{i}}{\weight(\Ipd) } \ip{z,z}, 
		\end{equation*}
		with $C_2 =  \frac12 \ps (\CML \e)^{\half}.$
	\end{lemma}
	\begin{proof}
	By \eqref{eq::Derivdh} we have $\abs{\grad \Ipd} \leq 1$. Thus, with definition \eqref{eq::DefWeightFunction} of $\weight(\cdot)$ we derive
	\begin{align*}
		\abs{\grad \weight(\Ipd) } = 
		\abs{\grad \exp\bigg( \frac{\gamma \Ipd}{\max \{ \lambda, h\}} \bigg)} 
		&= \abs{\frac{\gamma}{\max \{ \lambda, h\}} \exp\bigg( \frac{\gamma \Ipd}{\max \{ \lambda, h\}} \bigg) \grad \Ipd} \\
		&\leq  \frac{\gamma}{\max \{  \lambda , h\}} \weight(\Ipd) .
	\end{align*}
	The desired bound then follows by using again the \CSI, which yields
	\begin{align*}
		\abs{\int_{\SDov{i}} z \ps^2 \grad \weight(\Ipd) \cdot \grad z \dintx  } 
		&\leq \ps  \int_{\SDov{i}} \abs{z} \frac{\gamma}{\max \{  \lambda , h\}} \weight(\Ipd)  \ps \abs{\grad z} \dintx \\
		&\leq \frac{\ps \gamma }{\max \{  \lambda , h\}} \ipD{\SDov{i}}{\weight(\Ipd) z,z}^{\half}  \baw{\SDov{i}}{ \weight(\Ipd) }\ip{z,z}^{\half}\\
		&\leq \gamma \frac{\ps ( \CML \e^{\gamma})^{\half}}{2} \bbw{\SDov{i}}{\weight(\Ipd)} \ip{z,z} .
	\end{align*}
	The claim follows then as in \Cref{Lem::BoundDecayTerm2} from $\e^{\gamma} \leq \e$.
	\end{proof}

\section*{Acknowledgments}
We thank Pratik Kumbhar, Constantin Carle, and Benjamin D\"orich for many helpful discussions.

This work was funded by the Deutsche Forschungsgemeinschaft (DFG, German Research Foundation) — Project-ID 258734477 — CRC 1173. The authors acknowledge
support by the state of Baden-Württemberg through bwHPC.

\bibliographystyle{amsalpha}
\bibliography{references.bib}

\end{document}